\documentclass[3p]{elsarticle}
\journal{International Journal of Approximate Reasoning}

\usepackage{hyperref}
\usepackage{amssymb,amsfonts,amsmath,amsthm}

\newtheorem{theorem}{Theorem}[section]
\newtheorem{proposition}[theorem]{Proposition}
\newtheorem{lemma}[theorem]{Lemma}
\newtheorem{corollary}[theorem]{Corollary}

\newdefinition{definition}[theorem]{Definition}
\newdefinition{example}[theorem]{Example}
\newdefinition{remark}[theorem]{Remark}
\newdefinition{question}[theorem]{Question}

\newcommand{\UP}{\blacktriangle}
\newcommand{\DOWN}{\blacktriangledown}
\newcommand{\Up}{\vartriangle}
\newcommand{\Down}{\triangledown}

\newcommand{\Neg}{\daleth}

\begin{document}
\begin{frontmatter}

\title{Pseudo-Kleene algebras determined by rough sets}

\author{Jouni J{\"a}rvinen}
\address{Department of Software Engineering, {LUT}~University, Mukkulankatu~19, 15210~Lahti, Finland}
\ead{jouni.jarvinen@lut.fi}

\author{S{\'a}ndor Radeleczki}
\address{Institute of Mathematics, University of Miskolc, 3515~Miskolc-Egyetemv{\'a}ros, Hungary}
\ead{matradi@uni-miskolc.hu}

\begin{abstract}
 We study the pseudo-Kleene algebras of the Dedekind--MacNeille completion of
 the ordered set of rough set determined by a reflexive relation. We characterize 
 the cases when PBZ and PBZ*-lattices can be defined on these pseudo-Kleene algebras.
\end{abstract}

\begin{keyword}
  Rough set 
  \sep negation 
  \sep pseudo-Kleene algebra 
  \sep Brouwer--Zadeh lattice 
  \sep paraorhtomodular lattice 
  \sep Stone algebra 
\end{keyword}

\end{frontmatter}

\section{Introduction}
\label{Sec:Intro}

In rough set theory, introduced by Z.~Pawlak \cite{Pawl82},
knowledge about elements of a set $U$ is given in terms of an equivalence 
$E$ on $U$ interpreted so that $(x,y) \in E$ if the elements $x$ and $y$
cannot be distinguished in terms of the information represented by $E$. 
Each set $X \subseteq U$ is approximated by two sets: the lower approximation 
$X^\DOWN$ consists of elements which certainly belong to $X$ in view of knowledge $E$, 
and the upper approximation $X^\UP$ consists of objects 
which possibly are in $X$.

The pair $(X^\DOWN,X^\UP)$ is called a rough set. We denote by $\mathrm{RS}$ 
the set of all rough sets. It is proved in \cite{PomPom88} that $\mathrm{RS}$ is a Stone algebra.
This result was improved in \cite{Com93} by showing that $\mathrm{RS}$ forms a regular double 
Stone algebra. The three-valued Lukasiewicz algebras defined by RS were considered for the first time
in \cite{BanChak96}. P.~Pagliani showed in \cite{Pagliani97} how a 
semisimple Nelson algebra can be defined on $\mathrm{RS}$.

In the literature can be found studies in which the information about
the objects is given in terms of other types of relations than equivalences.
For instance, already in \cite{Yao96} rough approximations defined by an arbitrary binary relation $R$ on $U$
were considered. In that paper, the approximations were defined as in this study, that is,
for any $X \subseteq U$, an element $x$ belongs to the upper approximation $X^\UP$
whenever there is an element in $X$ to which $x$ is $R$-related. An element $x$ is in the lower approximation
$X^\DOWN$ if all elements to which $x$ is $R$-related are in $X$.
It is known that for an arbitrary tolerance (reflexive and symmetric binary relation), 
the ordered set $\mathrm{RS}$ is not necessarily a lattice; see \cite{Jarvinen2007}, for instance. 
In \cite{Umadevi2015}, D.~Umadevi presented the Dedekind--MacNeille completion of $\mathrm{RS}$ for arbitrary
binary relations. In this work, we denote this completion by $\mathrm{DM(RS)}$.

The work of Umadevi provides a starting point of for this study. Our aim is to find what
kind of logico-algebraic structures can be defined on $\mathrm{DM(RS)}$. We have
restricted ourselves to reflexive relations. Note that reflexivity is
equivalent to the fact that $X^\DOWN \subseteq X \subseteq X^\UP$ holds for any subset $X$ of $U$,   
which is a natural requirement for rough approximations.

Pseudo-Kleene algebras play an essential role in this study. They are bounded lattices equipped with a Kleene complement $\sim$. 
Note that  Kleene algebras are distributive pseudo-Kleene algebras. We show how to define a paraorthomodular pseudo-Kleene 
algebra on $\mathrm{DM(RS)}$. Interestingly, in the literature \cite{JarRad20,Qian2010} can be found studies in which approximation 
operators are defined by using certain combinations of two or more equivalence relations. In \cite{GegKovRad22}, the authors 
proved that the Dedekind--MacNeille completion of the so-called optimistic multigranular rough sets forms a
paraorthomodular pseudo-Kleene algebra.  Therefore, our work has a connection also to the study of multigranular approximation spaces.
This study also continues our lattice-theoretical research of $\mathrm{RS}$. In \cite{JarRad11} we proved that rough sets defined 
by quasiorders (reflexive and transitive binary relations) are exactly the Nelson algebras defined on algebraic lattices. 
In \cite{JarRad17} we characterized the rough sets defined by a tolerance  relation induced by an irredundant covering of $U$. 

In this paper, we show that the pseudo-Kleene algebras defined on $\mathrm{DM(RS)}$ have some distinctive properties. 
We prove that the sharp and the complemented elements coincide. 
We also note that if a complement of an element exists in $\mathrm{DM(RS)}$, it is unique. 
The pseudo-Kleene algebra $\mathrm{DM(RS)}$ is always paraorthomodular. 
Central elements have a key role in direct decompositions of bounded lattices, and
we establish a connection between central elements and exact rough sets.

Brouwer--Zadeh posets and lattices were introduced by G.~Cattaneo and  G.~Marino in \cite{Cattaneo84} in the setting of fuzzy sets. 
They were further investigated in \cite{Cattaneo89}, for example. Roughly speaking, they are structures with two complements: ${\sim}$
is the Kleene complementation and $\Neg$ behaves like an intuitionistic negation. Brouwer--Zadeh lattices
are commonly called BZ-lattices. Cattaneo studied BZ-lattices defined by preclusivity (irreflexive and symmetric) relations in \cite{Cattaneo97}. 
Together with D.~Ciucci he studied BZ-lattices related to rough sets structures determined by tolerances in \cite{CatCiu2005}. 

Paraorthomodular BZ-lattices are called PBZ-lattices. In this work, we study PBZ-lattices defined on $\mathrm{DM(RS)}$ in case of 
relations that are at least reflexive. In particular, we show that there is a one-to-one correspondence between atomistic 
complete Boolean sublattices of $\mathrm{DM(RS)}$ and PBZ-lattices on $\mathrm{DM(RS)}$.  
If $R$ a quasiorder or a tolerance induced by an irredundant covering, we
characterize the PBZ-lattices in terms of equivalences extending $R$.

PBZ*-lattices were introduced in \cite{Giuntini2016} in order to obtain insights into the structure of certain algebras of 
effects of a complex Hilbert space,  lattice-ordered by the so-called spectral order. They are 
PBZ-lattices satisfying the condition $\Neg(a \wedge {\sim} a) \leq 
\Neg a \vee \Neg {\sim} a$. We characterize the case when $\mathrm{RS}$ determined by a quasiorder defines a PBZ*-lattice. 

The paper is structured as follows. In Section~\ref{Sec:basic}, we
consider the basic properties of $\mathrm{DM(RS)}$.
Section~\ref{Sec:CentralElements} is devoted to 
central and exact elements of $\mathrm{DM(RS)}$. 
Basic properties of different kinds of Brouwer--Zadeh lattices are recalled 
from the literature in Section~\ref{sec:BZ-lattices}. Finally, in Section~\ref{sec:PBZ*-lattices}, 
we study PBZ and PBZ$^*$-lattices that can be defined on $\mathrm{DM(RS)}$. Some concluding remarks end the work.

\section{Smallest completion of rough sets}
\label{Sec:basic}

The rough sets lattices defined by equivalences, tolerances induced by irredundant coverings, and quasiorders 
are distributive, and they define a Kleene algebra. Pseudo-Kleene algebras
are non-distributive generalizations of Kleene algebras. 
We show that for a reflexive relation, a paraorthomodular pseudo-Kleene algebra can be defined on $\mathrm{DM(RS)}$. In this algebra, 
the sharp and the complemented elements coincide. Because this does not generally hold, 
we have that there are pseudo-Kleene algebras that are not isomorphic to some  $\mathrm{DM(RS)}$. 
We also show that if complements exist in $\mathrm{DM(RS)}$, 
they are unique. The section ends by presenting some lattice-theoretical
properties of the ordered set of the complemented elements. 

We begin by defining the rough set approximations based on arbitrary
binary relations. Let $R$ be a binary relation on $U$. We denote for
any $x \in U$, $R(x) := \{y \in U \mid (x,y) \in R \}$. 
The symbol $:=$ denotes `equals by definition'.
For any set $X \subseteq U$, the \emph{lower approximation} of $X$ is
\[ X^\DOWN := \{x \in U \mid R(x) \subseteq X \} \]
and the \emph{upper approximation} of $X$ is
\[ X^\UP := \{x \in U \mid R(x) \cap X \neq \emptyset\}.\]
Let $X^c$ denote the complement $U \setminus X$ of $X$. Then,
$X^{\DOWN c} = X^{c\UP}$ and $X^{\UP c} = X^{c \DOWN}$,
that is, $^\DOWN$ and $^\UP$ are dual.  

We may also determine rough set approximations in terms of the inverse
$R^{-1}$ of $R$, that is,
\[ X^\Down := \{x \in U \mid R^{-1}(x) \subseteq X \} \]
and
\[ X^\Up := \{x \in U \mid R^{-1}(x) \cap X \neq \emptyset \}. \]
Let $\wp(U)$ denote the family of all subsets of $U$.
Interestingly, the pairs $({^\UP},{^\Down})$ and $({^\Up},{^\DOWN})$ are 
order-preserving Galois connections on the complete lattice $(\wp(U),\subseteq)$. 
Several essential properties of the approximations follows from this
fact; see \cite{Jarvinen2007} for further details. In particular,
$\emptyset^\UP = \emptyset^\Up = \emptyset$ and 
$U^\DOWN = U^\Down = U$. 

In addition, we denote
\begin{gather} 
\wp(U)^\UP := \{ X^\UP \mid X \subseteq U\}, \quad
\wp(U)^\DOWN := \{ X^\DOWN \mid X \subseteq U\}, \\
\wp(U)^\Up := \{ X^\Up \mid X \subseteq U\}, \quad
\wp(U)^\Down := \{ X^\Down \mid X \subseteq U\}.
\end{gather}
For all $X \subseteq U$,
\begin{align*} 
X \in \wp(U)^\UP \iff X = X^{\Down\UP} \\
X \in \wp(U)^\DOWN \iff X = X^{\Up\DOWN}. 
\end{align*}

A binary relation $R$ on $U$ is \emph{right-total} if for any $x \in U$,
there is $y$ such that $(y,x) \in R$. This means that 
$R^{-1}(x) \neq \emptyset$ for all $x \in U$.
The notion of a \emph{left-total} is defined analogously, and obviously $R$ is left-total 
if and only if $R(x) \neq \emptyset$ for all $x \in U$. Note also
that left-total relations are often called \emph{serial}.

Let $R$ be right-total. Then clearly, $\emptyset^\Down = \emptyset$ and $U^\Up = U$.
Let $X \subseteq U$. If $X^{\Down \UP} = \emptyset$, then 
$X^\Down = X^{\Down \UP \Down} = \emptyset$. Similarly,
$X^{\Up \DOWN} = U$ implies $X^\Up = X^{\Up \DOWN \Up} = U^\Up = U$. 
For a left-total $R$, we have that  
$X^{\DOWN \Up} = \emptyset$ implies $X^{\DOWN} = \emptyset$
and $X^{\Up \DOWN} = U$ implies $X^{\Up} = U$.

\medskip\noindent%
Note that the following are equivalent:
\begin{enumerate}[(i)]
\item $R$ is reflexive;
\item $R^{-1}$ is reflexive;
\item $X \subseteq X^{\UP}$ and $X \subseteq X^{\Up}$ for all $X \subseteq U$;
\item $X^\DOWN \subseteq X$ and $X^\Down \subseteq X$ for all $X \subseteq U$.
\end{enumerate}

Let us denote by $\mathrm{RS}$ the set of all \emph{rough sets}, that is,
\[ \mathrm{RS} := \{ (X^\DOWN, X^\UP) \mid X \subseteq U \}. \]
The \emph{Dedekind--MacNeille completion} of an ordered set is the smallest complete lattice that contains it. 
We denote the Dedekind--MacNeille completion of $\mathrm{RS}$ by $\mathrm{DM(RS)}$.
Umadevi \cite{Umadevi2015} has proved that for any binary relation $R$ on $U$,
\[ \mathrm{DM(RS)} = \{ (A,B) \in \wp(U)^\DOWN \times \wp(U)^\UP \mid
A^{\Up \UP} \subseteq B \ \text{ and } \ A \cap \mathcal{S} 
= B \cap \mathcal{S} \}. \]
Here $\mathcal{S}$ is the set of \emph{singletons}, that is,
\[ 
\mathcal{S} := \{ x \in U \mid R(x) = \{z\} \text{ for some $z \in U$} \}. 
\]
For $\{(A_i,B_i)\}_{i \in I} \subseteq \mathrm{DM(RS)}$,
\begin{align*}
\bigvee_{i \in I} (A_i,B_i) 
&= \Big ( \big ( \bigcup_{i \in I} A_i \big )^{\Up \DOWN}, 
\bigcup_{i \in I} B_i \Big )\\ 
\bigwedge_{i \in I} (A_i,B_i)  & = 
\Big ( \bigcap_{i \in I}A_i, \big (\bigcap_{i \in I} B_i \big )^{\Down \UP} \Big). 
\end{align*}
Umadevi also showed that for $(A,B) \in \mathrm{DM(RS)}$, the pair $(B^c, A^c)$
belongs to $\mathrm{DM(RS)}$. This means that the map
\[ {\sim} \colon  \mathrm{DM(RS)} \to  \mathrm{DM(RS)}, (A,B) \mapsto (B^c, A^c) \]
forms an antitone involution, that is,
\begin{gather}
{\sim}{\sim}(A,B) = (A,B) \\
(A,B) \leq (C,D) \text{ \ implies \ } {\sim}(C,D) \leq {\sim}(A,B)
\end{gather}

\medskip

In this work, we consider rough approximations defined by a reflexive relation. This
means that for any $(A,B) \in \mathrm{DM(RS)}$,
$A \subseteq A^{\Up \UP} \subseteq B$.

A \emph{pseudo-Kleene poset} $(P,\leq,{\sim})$ is an ordered set $(P,\leq)$ 
equipped with an antitone involution $\sim$ satisfying the condition:
\begin{equation} \label{Eq:psKleenePoset}
\text{if $a \leq {\sim} a$  and  ${\sim} b \leq b$, then $a \leq b$.}    
\end{equation}
for all $a,b \in P$.
Note that if $(L,\leq)$ forms a lattice, then \eqref{Eq:psKleenePoset}
is equivalent to condition
\begin{equation} \label{Eq:psKleeneAlg}
a \wedge {\sim}a \leq b \vee {\sim} b.     
\end{equation}
The algebra $(L,\vee,\wedge,{\sim})$ is called a \emph{pseudo-Kleene algebra}.

Umadevi  \cite{Umadevi2015} mentioned without proof that for any binary relation,
$\mathrm{DM(RS)}$ forms a pseudo-Kleene algebra. For the sake of completeness,
we write the following lemma.

\begin{lemma} \label{lem:pseudo_kleene}
Let $R$ be reflexive binary relation on $U$. Then
\begin{enumerate}[\rm (a)]
\item $\mathrm{RS}$ forms a pseudo-Kleene poset;
\item $\mathrm{DM(RS)}$ forms a pseudo-Kleene algebra.
\end{enumerate}
\end{lemma}

\begin{proof} Let $(A,B), (C,D) \in \mathrm{DM(RS)}$. Since $R$ is reflexive, 
$A \subseteq B$ and $C \subseteq D$. We have $B^c \subseteq A^c$ and
$D^c \subseteq C^c$. These give $A \cap B^c \subseteq A \cap A^c = \emptyset$
and $U = D \cup D^c \subseteq D \cup C^c$.

(a) If $(A,B) \leq {\sim}(A,B)$, then ${\sim}(A,B) = (B^c, A^c)$ implies 
$A \subseteq B^c$. We obtain $A \subseteq B \cap B^c = \emptyset$. 
If ${\sim}(C,D) \leq (C,D)$, then ${\sim}(C,D) = (D^c,C^c)$ gives
$C^c \subseteq D$ and  $U = C^c \cup C \subseteq D$. We have
\[ (A,B) = (\emptyset,B) \leq (C,U) = (C,D).\]

(b) Similarly,
\begin{gather*}
    (A,B) \wedge {\sim} (A,B)  = (A,B) \wedge (B^c,A^c) =
    (A \cap B^c, (B \cap A^c)^{\Down \UP})
    = (\emptyset, (B \cap A^c)^{\Down \UP}) \\ 
     \leq ((C \cup D^c)^{\Up \DOWN},U) 
     = ((C \cup D^c)^{\Up \DOWN}, D \cup C^c) 
     = (C,D) \vee (D^c,C^c) = (C,D) \vee {\sim} (C,D).   \qedhere
\end{gather*}   
\end{proof} 

A bounded pseudo-Kleene algebra $(L,\vee,\wedge, {\sim},0,1)$ is said to be
\emph{paraorthomodular} if for all $a,b \in L$,
\[ a \leq b \ \text{and} \ {\sim}a \wedge b = 0 \ \text{imply} \ a = b.\]

\begin{proposition} \label{Prop:RS_paraorthomodular}
If $R$ is reflexive, then $\mathrm{DM(RS)}$ is paraorthomodular.
\end{proposition}

\begin{proof}
Suppose that $(A,B) \leq (C,D)$ and 
${\sim}(A,B) \wedge (C,D) = (\emptyset,\emptyset)$ in  $\mathrm{DM(RS)}$.
The first condition means $A \subseteq C$ and $B \subseteq D$. 
From the second condition, we get 
\[
    (B^c, A^c) \wedge (C,D) = (B^c \cap C, (A^c \cap D)^{\Down \UP}) = 
    (\emptyset,\emptyset).
\]
Now $B^c \cap C = \emptyset$ implies 
$C \subseteq B^{cc} = B$. In addition, $(A^c \cap D)^{\Down \UP} = \emptyset$
gives $\emptyset = (A^c \cap D)^{\Down}  = 
A^{c \Down} \cap D^\Down = 
A^{\Up c} \cap D^\Down$.
This is equivalent to $D^\Down \subseteq A^{\Up c c} = A^\Up$. Because
$D \in \wp(U)^\UP$, 
\[ D = D^{\Down \UP} \subseteq A^{\Up \UP} \subseteq B.\]
As we already noted, $B \subseteq D$. Therefore, $B = D$.

Since $(C,D) \in \mathrm{DM(RS)}$, $C^{\Up \UP} \subseteq D$. This means that
for all $x \in U$, $R(x) \cap C^\Up \neq \emptyset$ implies $x \in D$.
Let $y \in C^\Up$ and $z \in R^{-1}(y)$. Because $y \in R(z) \cap C^\Up$, we have $z \in D$. 
Thus, $R^{-1}(y) \subseteq D$ and $y \in D^\Down$.
We have now shown that  $C^\Up \subseteq D^\Down$.
We already noted that $D^\Down \subseteq A^\Up$. Thus, 
$C^\Up \subseteq A^\Up$. The fact that $A \subseteq C$ gives 
$A^\Up \subseteq C^\Up$. Thus, $A^\Up = C^\Up$. Because $A$ and $C$
belong to $\wp(U)^\DOWN$, we obtain $A = A^{\Up \DOWN} = C^{\Up \DOWN} = C$.

We have now proved that $(A,B) = (C,D)$.  
\end{proof}

\begin{remark} \label{rem:Chajda}
In \cite{Chajda2016}, I.~Chajda defined pseudo-Kleene algebras as
lattices with an antitone involution $\sim$ satisfying
\eqref{Eq:psKleeneAlg} and
\begin{equation} \label{eq:semi_distr}
x \wedge ({\sim}x \vee y) =
(x \wedge {\sim}x) \vee (x \wedge y).
\end{equation}

Let $R$ be a reflexive relation on $U = \{a,b,c\}$ such that
$R(a) = \{a,b\}$, $R(b) = \{b,c\}$, and $R(c) = \{c\}$.
The lattice $\mathrm{RS}$ is depicted in Figure~\ref{fig:twoFigs}(a).
Note that in the figure, sets are denoted by sequences by their elements. For
instance, $\{a,b\}$ is denoted by $ab$.
\begin{figure} 
\centering
\begin{minipage}{.5\textwidth}
  \centering
  \includegraphics[height=55mm]{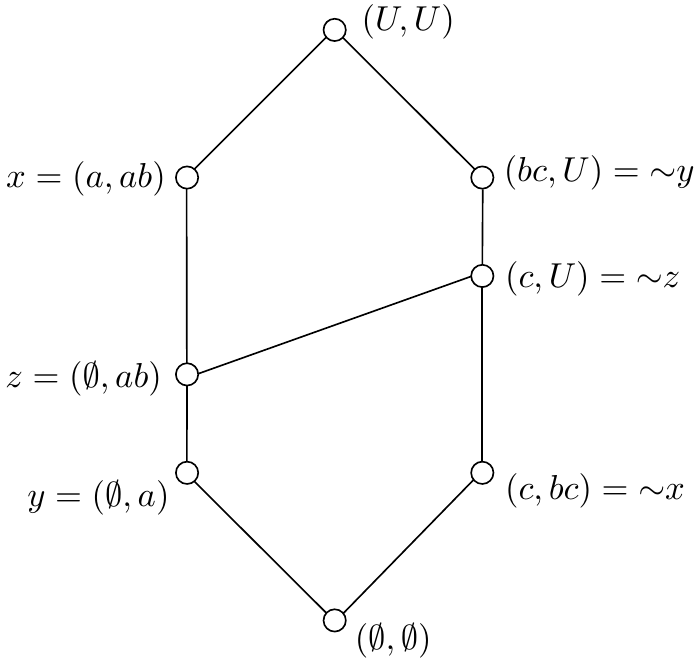} 
\end{minipage}%
\begin{minipage}{.5\textwidth}
  \centering
  \includegraphics[height=55mm]{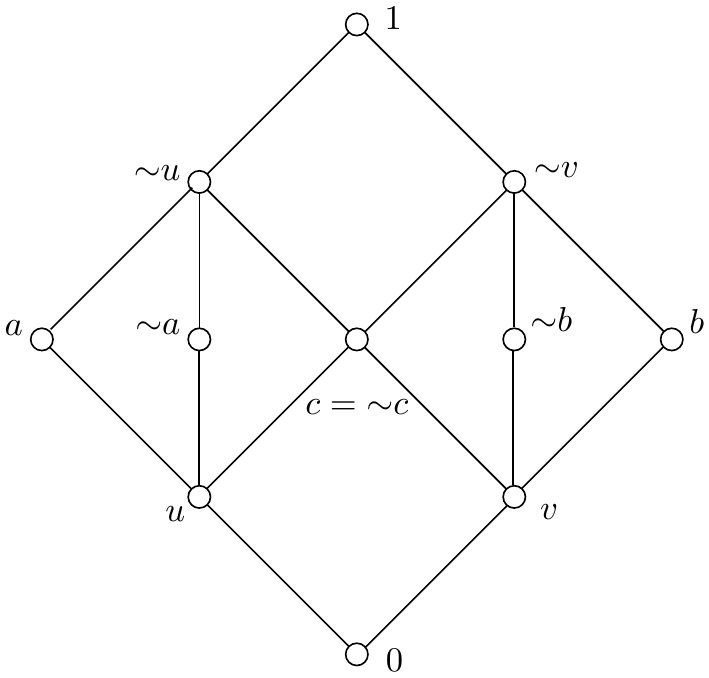}
\end{minipage}
 \caption{(a) Lattice $\mathrm{RS}$ of Remark~\ref{rem:Chajda}. (b) Pseudo-Kleene algebra of Example~\ref{exa:muresan}. \label{fig:twoFigs}}
\end{figure}

Let us set
\[ x := (\{a\}, \{a,b\}), \quad
   z:= (\emptyset,\{a,b\}), \quad
   y:= (\emptyset,\{a\}). \]
Now 
\[ x \wedge ({\sim}x \vee y) = z \text{ and }
  (x \wedge {\sim}x) \vee (x \wedge y) = y 
\]
This means that \eqref{eq:semi_distr} does not generally hold in
$\mathrm{RS}$ or $\mathrm{DM(RS)}$.
\end{remark}

In lattice-theory, an element $a$ of a bounded lattice $L$ is called \emph{complemented} if there 
is an element $b \in L$ such that 
\[ a \wedge b = 0 \text{ \ and \ } a \vee b = 1.\]
An element $a$ of a bounded pseudo-Kleene algebra is called \emph{sharp} if $a \wedge {\sim}a = 0$. 
It is easy to see that every sharp element is complemented. Indeed, let $a$ be a sharp element. 
Then $a \wedge {\sim}a = 0$ and this implies
\[ a \vee {\sim}a = {\sim}{\sim}a \vee {\sim}a = 
{\sim}({\sim}a \wedge a) = {\sim}0 = 1.\]
As noted in \cite{Giuntini2016}, if $a \leq b$ and ${\sim} a \wedge b = 0$,
then $a$ and $b$ are sharp. Indeed, $0 = a \wedge 0 = 
a \wedge ({\sim}a \wedge b) = (a \wedge b) \wedge {\sim}a =
a \wedge {\sim}a$. Note that ${\sim} a \wedge b = 0$ is equivalent to
$a \vee {\sim} b = 1$. Therefore, $1 = b \vee 1 = b \vee (a \vee {\sim}b) =
(a \vee b) \vee {\sim}b = b \vee {\sim}b$.

\begin{example} \label{exa:muresan}
Not every complemented element is sharp, as can be seen in the pseudo-Kleene algebra 
of Figure~\ref{fig:twoFigs}(b), which appears originally in 
\cite[Example 3.5]{Muresan2019}.

Now the element $a$ is complemented, it has two complements $b$ and
${\sim}b$. However, $a$ is not sharp, because $a \wedge {\sim}a = u \neq 0$.
\end{example}

Example~\ref{exa:muresan} reveals that there are complemented elements which are not sharp.  It should be noted that already in \cite{Umadevi2015}, 
Umadevi characterized the complemented elements. Therefore, our next 
proposition can be seen as an extension of her result.

\begin{proposition} \label{Prop:SharpCharacterization}
Let $R$ be a reflexive relation on $U$. For $(A,B)$ in
$\mathrm{DM(RS)}$, the following are equivalent:
\begin{enumerate}[\rm (a)]
    \item $B^\Down = A^\Up$;
    \item $(A,B)$ is sharp;
    \item $(A,B)$ is complemented.
\end{enumerate}
\end{proposition}

\begin{proof}
(a)$\Rightarrow$(b): Suppose $B^\Down = A^\Up$. Now
\[ (A,B) \wedge {\sim}(A,B) = (A,B) \wedge (B^c, A^c) = (A \cap B^c, (B \cap A^c)^{\Down\UP}).\]
Because $R$ is reflexive, $A \subseteq B$ and $B^c \subseteq A^c$. Thus,
$A \cap B^c \subseteq A \cap A^c = \emptyset$. 

Secondly, $B^\Down \subseteq A^\Up$ implies
$A^{\Up c} \subseteq B^{\Down c}$ and hence
\[    (B \cap A^c)^{\Down}
= B^\Down \cap A^{c \Down} = B^\Down \cap A^{\Up c} 
\subseteq B^\Down \cap B^{\Down c} = \emptyset. \]
Since $(B \cap A^c)^{\Down} = \emptyset$, we have
$(B \cap A^c)^{\Down\UP} = \emptyset^\UP = \emptyset$.

Thus, $ (A,B) \wedge {\sim}(A,B) = (\emptyset, \emptyset)$ and  $(A,B)$ is sharp.

\medskip\noindent%
(b)$\Rightarrow$(c): We have already noted that in a bounded pseudo-Kleene algebra, each sharp element is complemented.

\noindent%
(c)$\Rightarrow$(a): Let $(A,B) \in \mathrm{DM(RS)}$. Then $A^{\Up \UP} \subseteq B$, which implies $A^\Up \subseteq (A^{\Up})^{\UP \Down} \subseteq B^\Down$. 
Suppose $(A,B)$ is complemented. Then there exists $(C,D)$ in 
$\mathrm{DM(RS)}$ such that
\begin{align}
(A,B) \wedge (C,D) = (\emptyset,\emptyset) \label{eq:wedge}\\ 
\intertext{and} 
(A,B) \vee (C,D) = (U,U) . \label{eq:vee}
\end{align}
We need to prove that $B^\Down \subseteq A^\Up$. From
\eqref{eq:vee} we get $(A\cup C)^{\Up \DOWN} = U$, which is 
equivalent to $U = (A\cup C)^{\Up} = A^\Up \cup C^\Up$. This means 
$ A^{c \Down} = A^{\Up c} \subseteq C^\Up$.
Because $A \in \wp(U)^\DOWN$, $A = A^{\Up \DOWN}$. We obtain
$A^c = A^{\Up \DOWN c} = A^{c \Down \UP} \subseteq C^{\Up \UP} \subseteq D$.

From \eqref{eq:wedge} we have that $(B \cap D)^{\Down \UP} = \emptyset$. 
This is equivalent to $\emptyset = (B \cap D)^\Down = B^\Down \cap D^\Down$,
that is, $D^\Down \subseteq B^{\Down c}$. We can now write
\[ A^{c \Down} \subseteq D^\Down \subseteq B^{\Down c} \]
and
\[ B^\Down \subseteq A^{c \Down c} = A^{\Up c c} = A^\Up,\]
which completes the proof.  
\end{proof}

Let $\mathcal{C}$ denote the set of complemented elements of DM$(RS)$. 
By Proposition~\ref{Prop:SharpCharacterization}, $\mathcal{C}$ is also the
set of the sharp elements. 

\begin{lemma} \label{lem:simple_sharp}
Let $R$ be a reflexive relation on $U$.
\begin{enumerate}[\rm (a)]
\item $\mathcal{C} \subseteq \mathrm{RS}$.
\item $(\emptyset,\emptyset)$, $(U,U) \in \mathcal{C}$.
\item If $(A,B) \in \mathcal{C}$, then ${\sim}(A,B) \in \mathcal{C}$.
\end{enumerate}
\end{lemma}

\begin{proof}
(a) Suppose that $(A,B) \in \mathrm{DM(RS)}$ is complemented. 
Then, by Proposition~\ref{Prop:SharpCharacterization},
$A^\Up = B^\Down$. Let us denote $X := A^\Up = B^\Down$. Now
$X^\DOWN = A^{\Up \DOWN} = A$ and $X^{\UP} = B^{\Down \UP} = B$. Thus,
$(A,B) = (X^\DOWN, X^\UP)$.

(b) Since $R$ is reflexive, it is well-known and obvious that
$(\emptyset,\emptyset) \in \mathrm{RS}$ and $\emptyset^\Up = \emptyset^\Down = \emptyset$. 
Similarly, $(U,U) \in \mathrm{RS}$ and $U^\Up = U^\Down = U$.

(c) For any $(A,B) \in \mathcal{C}$, we
have $B^{c \Up} = B^{\Down c} = A^{\Up c} = A^{c \Down}$, meaning
that ${\sim}(A,B) = (B^c, A^c)$ is in $\mathcal{C}$.
\end{proof}

Our following proposition shows that if an element has a complement in
$\mathrm{DM(RS)}$, it is unique.

\begin{proposition}\label{Prop:UniqueComplement}
Let $R$ be a reflexive relation on $U$ and $(A,B) \in \mathrm{DM(RS)}$. 
If $(A,B)$ has a complement $(C,D)$, then $(C,D) = {\sim}(A,B)$.
\end{proposition}

\begin{proof} Suppose that $(C,D)$ is a complement of $(A,B)$. As in
the proof of Proposition~\ref{Prop:SharpCharacterization}, we have
$A^c \subseteq D$, $A^\Up = B^\Down$ and $A^\Up \cup C^\Up = U$. Now
$B^{c \Up} = B^{\Down c} = A^{\Up c} \subseteq C^\Up$. Since $C \in \wp(U)^\DOWN$,
we have $C = C^{\Up \DOWN}$. This yields
$B^c \subseteq B^{c \Up \DOWN} \subseteq C^{\Up \DOWN} = C$.

We have now shown ${\sim}(A,B) = (B^c,A^c) \leq (C,D)$. 
Because $\mathrm{DM(RS)}$ is a paraorthomodular
lattice with $(\emptyset,\emptyset)$ as the least element, 
${\sim}(A,B) \leq (C,D)$ and ${\sim} {\sim}(A,B) \wedge (C,D) = (\emptyset,\emptyset)$
imply $(C,D) = {\sim}(A,B)$.
\end{proof}

Let $L$ be a lattice in which the complementation $'$ is unique.
It is known \cite[Theorem 6.5]{Blyth} that if  $x \leq y$ implies $y' \leq x'$
for all $x,y \in L$, then $L$ is distributive. By this fact, we can write the following corollary.

\begin{corollary} \label{cor:BoolComplement}
If $\mathcal{C}$ is a sublattice of $\mathrm{DM(RS)}$, then $\mathcal{C}$ is Boolean.    
\end{corollary}

Our next observation shows that $\mathcal{C}$ is not always a sublattice of
$\mathrm{DM(RS)}$.

\begin{example}\label{exa:not_sublattice}
Let $U = \{1,2,3,4,5\}$ and $R$ a tolerance on $U$ such that
$R(1) = \{1,2\}$, $R(2) = \{1,2,3\}$, $R(3) = \{2,3,4\}$,
$R(4) = \{3,4,5\}$, and $R(5) = \{4,5\}$. It is known \cite{Jarv99} that the 
ordered set $\mathrm{RS}$ defined by $R$ is not a lattice.

The completion $\mathrm{DM(RS)}$ is given in Figure~\ref{fig:completed}.
The two elements which belong to $\mathrm{DM(RS)}$ but not in
$\mathrm{RS}$ are inside small boxes.
 
\begin{figure}[ht]
\centering
\includegraphics[width=150mm]{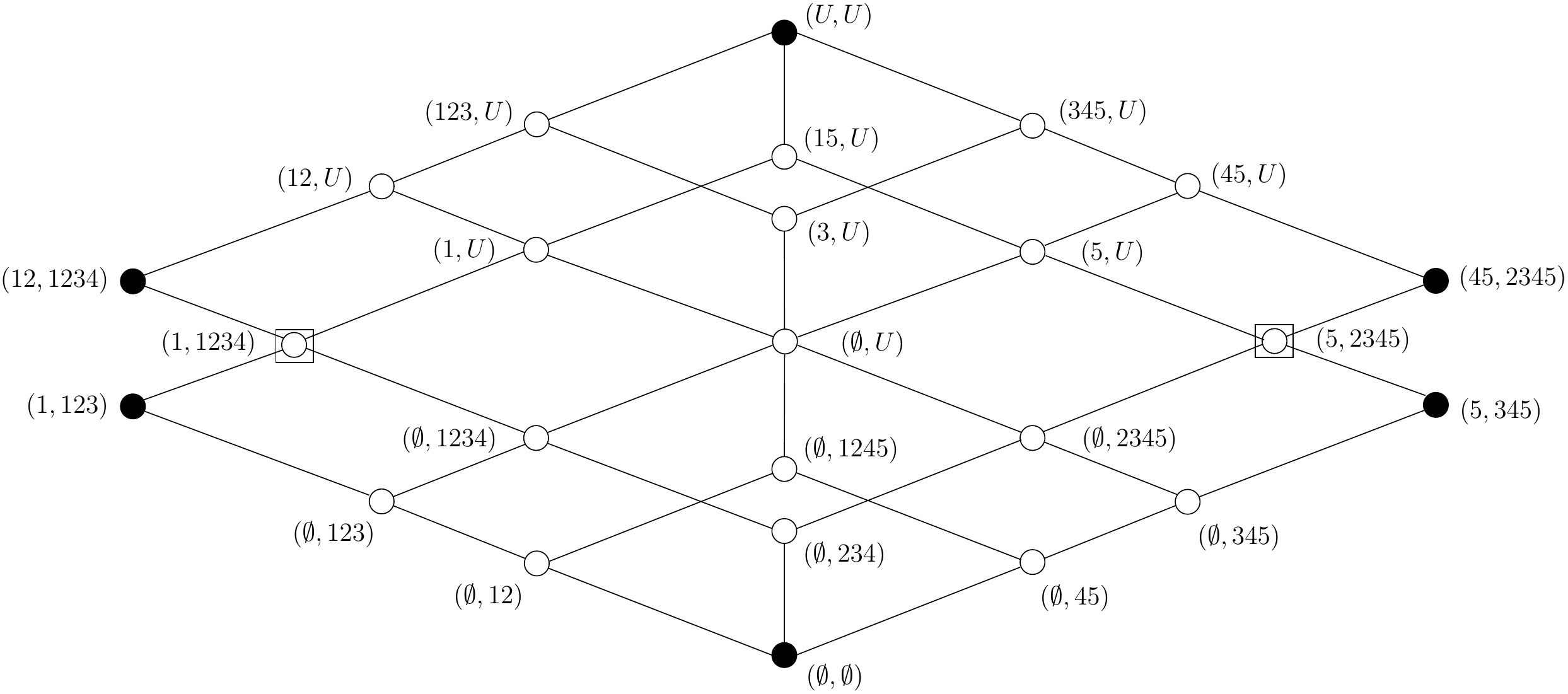}
\caption{The lattice  $\mathrm{DM(RS)}$ \label{fig:completed}}
\end{figure}

 The complemented elements are marked with filled circles. Clearly, $\mathcal{C}$ is not
 a sublattice of $\mathrm{DM(RS)}$. Also $\mathcal{C}$ is not distributive,
 because it contains $\mathbf{N_5}$ as a sublattice. One can also observe that
 as a lattice, $\mathcal{C}$ is not uniquely complemented.
\end{example}

Lemma~\ref{lem:simple_sharp} implies that by restricting the partial order $\leq$ of $\mathrm{DM(RS)}$ 
to $\mathcal{C}$, we obtain a bounded poset with involution $(\mathcal{C},\leq , {\sim})$. 
Let us consider the set
\[
\mathcal{A} := \wp(U)^\Down \cap \wp(U) ^\Up.
\]

\begin{lemma}\label{Lemma:A-set}
Let $R$ be a reflexive relation on $U$.
\begin{enumerate}[\rm (i)]
\item If $Z \in \mathcal{A}$, then $Z^c \in \mathcal{A}$.
\item $\mathcal{A} = \{Z \subseteq U\mid Z^{\DOWN
\Up}=Z^{\UP\Down} \}$.
\item  $(A,B)\in \mathrm{DM(RS)}$ is sharp if and only if $(A,B)=(Z^\UP,Z^\DOWN)$ for some $Z \in \mathcal{A}$.
\end{enumerate}
\end{lemma}

\begin{proof}
(i) Let $Z \in \mathcal{A}$. Then there are $X,Y \subseteq U$ such that $Z = X^\Up = Y^\Down$. 
Now, $Z^{c} = X^{\Up c} = ( X^{c})^\Down$ and $Z^{c} = Y^{\Down c} = (Y^{c})^\Up$ imply that $Z^c \in \wp(U)^{\Down} \cap \wp(U) ^{\Up} = \mathcal{A}$.

(ii) For an arbitrary binary relation $R$, the inclusions $Z^{\DOWN \Up }\subseteq Z\subseteq Z^{\UP \Down}$ hold. 
Hence, $Z^{\DOWN \Up} = Z^{\UP \Down}$ implies $Z = Z^{\DOWN \Up} = Z^{\UP \Down}$, that is,
$Z \in \mathcal{A}$. We have shown that 
$\{Z\subseteq U\mid  Z^{\DOWN \Up} = Z^{\UP \Down} \}\subseteq \mathcal{A}$.

Conversely, let $Z \in \mathcal{A}$.
Then, $Z = X^{\Up} = Y^{\Down}$ for some $X,Y\subseteq U$. 
We obtain $Z^{\DOWN \Up}=X^{\Up \DOWN \Up} = X^{\Up} = 
Y^{\Down }=Y^{\Down \UP \Down} = Z^{\UP \Down}$, proving that
$\mathcal{A} \subseteq \{Z\subseteq U\mid  Z^{\DOWN \Up} = Z^{\UP \Down} \}$.

(iii) Assume that $(A,B)\in  \mathrm{DM(RS)}$ is sharp. Let
$Z := A^{\Up} = B^{\Down}$. Then, $Z\in \wp(U)^{\Down} \cap \wp(U)^{\Up} = \mathcal{A}$. Now $Z^\DOWN = A^{\Up \DOWN} = A$ and
$Z^\UP = B^{\Down \UP} = B$, that is, $(A,B)=(Z^\DOWN,Z^\UP)$. 

On the other hand, suppose that $(A,B) = (Z^\DOWN, Z^\UP)$ for some 
$Z\in \mathcal{A}$. Then, $(A,B) \in \mathrm{RS} \subseteq \mathrm{DM(RS)}$. Now $A^{\Up} = Z^{\DOWN \Up} = Z^{\UP \Down} = B^{\Down }$. 
This means that $(A,B)$
is sharp.
\end{proof}

\begin{corollary}\label{Cor: A}
If $R$ is a reflexive relation on $U$, then $(\mathcal{A},{\subseteq}, {^c})$ is a bounded pseudo-Kleene poset.
\end{corollary}

\begin{proof}
It is clear that $\emptyset ,U\in \mathcal{A}$. Hence,  $(\mathcal{A},\subseteq)$ is a bounded. By Lemma~\ref{Lemma:A-set},
$\mathcal{A}$ is closed under $^c$. 
Let $Z_1,Z_2 \in \mathcal{A}$ be such that $Z_1 \subseteq {Z_1}^c$ and 
${Z_2}^c \subseteq Z_2$. Then $Z_1 \subseteq Z_1 \cap {Z_1}^c = \emptyset$ and 
$U = Z_2 \cup {Z_2}^c \subseteq Z_2$ give $Z_1 \subseteq Z_2$.
\end{proof}

It is also clear that for $Z \subseteq U$, $Z \in \mathcal{A}$
if and only if $Z = Z^{\DOWN \Up}$ and $Z = Z^{\UP \Down}$.
 
\begin{proposition} \label{Prop:isomorphism}
If $R$ is a reflexive relation on $U$, then 
$(\mathcal{C},{\leq} ,{\sim} )$ and $(\mathcal{A},{\subseteq}, {^c})$ are isomorphic pseudo-Kleene posets.
\end{proposition}

\begin{proof}
Let us define the map $\varphi \colon \mathcal{C} \rightarrow \mathcal{A}$
by setting $\varphi (A,B) = A^\Up $ for all $(A,B)\in \mathcal{C}$. 
As $A^\Up = B^\Down$ for any $(A,B) \in \mathcal{C}$, we have $\varphi (A,B) \in \mathcal{A}$. Thus, the map $\varphi$ is well-defined.

Let $(A_{1},B_{1}),(A_{2},B_{2}) \in \mathcal{C}$.
If $(A_{1},B_{1}) \leq (A_{2},B_{2})$, then $A_{1}\subseteq A_{2}$,
whence we get $\varphi (A_{1},B_{1}) = {A_1}^\Up
\subseteq {A_2}^\Up = \varphi (A_{2},B_{2})$.

On the other hand, if $\varphi (A_{1},B_{1}) = {A_1}^\Up
\subseteq {A_2}^\Up = \varphi (A_{2},B_{2})$, then 
$A_1 = {A_1}^{\Up \DOWN} \subseteq {A_2}^{\Up \DOWN} = A_2$ and
$B_1 = {B_1}^{\Down \UP} = {A_1}^{\Up \UP} \subseteq {A_2}^{\Up \UP}
= {B_2}^{\Down \UP} = B_2$. This means
that $(A_{1},B_{1}) \leq (A_{2},B_{2})$. We have shown that $\varphi$ is
an order-embedding.

Suppose that $Z \in \mathcal{A}$. By Lemma~\ref{Lemma:A-set},
$(Z^\DOWN,Z^\UP)$ belongs to $\mathcal{C}$. We have
$\varphi(Z^\DOWN,Z^\UP) = Z^{\DOWN \Up} = Z$. This means that 
$\varphi$ is onto $\mathcal{A}$. We have now proved that 
$\varphi$ is an order-isomorphism.

Let $(A,B) \in \mathcal{C}$. Then $B^\Down = A^\Up$ and 
\[ 
\varphi ({\sim} (A,B)) = \varphi  (B^{c},A^{c} ) = B^{c\Up} = B^{\Down c} = A^{\Up c} = \varphi(A,B)^{c}. \qedhere
\]
\end{proof}

Let us define a map $\psi \colon \mathcal{A} \to \mathcal{C}$ by
$\psi(Z) = (Z^\DOWN,Z^\UP)$ for $Z \in \mathcal{A}$. By Lemma~\ref{Lemma:A-set}(iii), $\Psi(Z) \in \mathcal{C}$, so the
map is well-defined. It is also clear that $\psi$ is order-preserving.
Indeed, if $Z_1, Z_2 \in \mathcal{A}$ are such that $Z_1 \subseteq Z_2$,
then $\psi(Z_1) = ({Z_1}^\DOWN,{Z_1}^\UP) \leq ({Z_2}^\DOWN,{Z_2}^\UP) =
\psi(Z_2)$.

\begin{lemma} Let $R$ be a reflexive relation on $U$. 
The map $\psi$ is the inverse of $\varphi$.
\end{lemma} 

\begin{proof} Let $(A,B) \in \mathcal{C}$. Because $A^\Up = B^\Down$,
we have
\[
\psi(\varphi(A,B)) = \psi(A^\Up) = (A^{\Up \DOWN}, A^{\Up \UP}) 
=  (A^{\Up \DOWN}, B^{\Down \UP}) = (A,B) 
\]
Note that the last equality follows from the facts that $A \in \wp(U)^\DOWN$ and $B \in \wp(U)^\UP$.

On the other hand, if $Z \in \mathcal{A}$, then
\[ \varphi ( \psi (Z)) = \varphi(Z^\DOWN, Z^\UP) = 
Z^{\DOWN \Up} = Z. \qedhere \]
\end{proof}

Suppose that $\psi(Z_1) \leq \psi(Z_2)$. Because $\varphi$ is order-preserving, 
$Z_1 = \varphi(\psi(Z_1)) \leq \varphi(\psi(Z_2)) = Z_2$. As
we already noted, $\psi$ is order-preserving. Thus, $\psi$ is an
order-embedding. It is clear that $\psi$ is a bijection, so it is
an order-isomorphism.

\begin{remark}\label{Rem:sublattice}
Let $\mathcal{T} \subseteq \mathcal{C}$ be closed under $\sim$. 
Then, $\mathcal{T}$ is a complete sublattice of $\mathrm{DM(RS)}$ if and only if 
it is closed with respect to arbitrary joins. 

Indeed, assume that $\mathcal{T}$ is closed with respect to arbitrary
joins and let $\mathcal{H} \subseteq  \mathcal{T}$. 
Then ${\sim} Z \in \mathcal{T}$ for each $Z \in \mathcal{H}$,
and hence $\bigvee \{ {\sim} Z \mid Z \in \mathcal{H} \}$ belongs to $\mathcal{T}$. We obtain 
\[ \bigwedge \mathcal{H}  =
{\sim} \bigvee \{ {\sim} Z \mid Z \in \mathcal{H} \}
\in \mathcal{T},
\]
proving that $\mathcal{T}$ is a complete sublattice of 
$\mathrm{DM(RS)}$. 

Analogously, let $\mathcal{S} \subseteq \mathcal{A}$ be closed under
set-complementation. Then $\mathcal{S}$ is a complete sublattice of $\wp(U)$ 
if and only if it is closed with respect to arbitrary unions. 
\end{remark}

For a function $f \colon A \to B$,  we write for all $S \subseteq A$,
\[ f[S] := \{f(x) \mid x \in S \} .\] 

\begin{proposition} \label{Prop:Sharplattice}
Let $R$ be a reflexive relation on $U$ and 
$\mathcal{T} \subseteq \mathcal{C}$ a nonempty set closed with respect to $ \sim $. 
Then $\mathcal{T}$ is a complete sublattice of $\mathrm{DM(RS)}$ if and only if $\varphi[\mathcal{T}]$ is a complete sublattice of 
$\wp(U)$.
\end{proposition}

\begin{proof}
Assume that $\mathcal{T}$ is a complete sublattice of 
$\mathrm{DM(RS)}$ and let $\mathcal{H} \subseteq \varphi[\mathcal{T}]$. 
Then $\psi(Z) = (Z^\DOWN,Z^\UP)$ for all $Z \in \mathcal{H}$.
Since each $\psi(Z)$ belongs to $\mathcal{T}$ and $\mathcal{T}$ is a complete sublattice of $\mathrm{DM(RS)}$, we have
\[ \bigvee_\mathcal{T} \{ \psi(Z) \mid Z \in \mathcal{H} \}
=  \bigvee_\mathcal{T} \{  (Z^\DOWN, Z^\UP) \mid Z \in \mathcal{H} \}
=  \Big ( \big ( \bigcup_{Z \in \mathcal{H}} Z^\DOWN \big)^{\Up \DOWN },
\bigcup_{Z \in \mathcal{H}} Z^\UP \Big) \in \mathcal{T}.
\]
Now $\varphi[\mathcal{T}]$ is a complete lattice isomorphic to $\mathcal{T}$ and
\[ \varphi \Big ( 
\bigvee_\mathcal{T} \{ \psi(Z) \mid Z \in \mathcal{H} \}
\Big )
= \bigvee_{\varphi[\mathcal{T}]} \{ \varphi(\psi(Z)) \mid Z \in \mathcal{H} \}
= \bigvee_{\varphi[\mathcal{T}]} \{ Z \mid Z \in \mathcal{H} \}
= \bigvee_{\varphi[\mathcal{T}]} \mathcal{H}.
\]
This implies that
\[ \bigvee_{\varphi[\mathcal{T}]} \mathcal{H} = 
\varphi \Big ( \big ( \bigcup_{Z \in \mathcal{H}} Z^\DOWN \big)^{\Up \DOWN },
\bigcup_{Z \in \mathcal{H}} Z^\UP \Big)
= \big ( \bigcup_{Z \in \mathcal{H}} Z^\DOWN \big)^{\Up \DOWN \Up}
= \big ( \bigcup_{Z \in \mathcal{H}} Z^\DOWN \big)^{\Up}
= \bigcup_{Z \in \mathcal{H}} Z^{\DOWN \Up}
= \bigcup_{Z \in \mathcal{H}} Z.
\]
In view of Remark~\ref{Rem:sublattice} this means that $\varphi[\mathcal{T}]$ is a complete sublattice of $\wp(U)$.

Conversely, assume that $\varphi[\mathcal{T}]$ is a complete sublattice of $\wp(U)$. Let $\mathcal{H} \subseteq \mathcal{T}$. 
Then $A^\Up = B^\Down$ for all $(A,B) \in \mathcal{H}$. Let us consider the join $\bigvee \mathcal{H}$
formed in $\mathrm{DM(RS)}$. We have 
\[ \bigvee \mathcal{H} = (\alpha, \beta), \]
where
\[ \alpha = \big ( \bigcup { A \mid (A,B) \in \mathcal{H}} \} \big )^{\Up \DOWN}
\qquad \text{and} \qquad
\beta = \bigcup \{B \mid (A,B) \in \mathcal{H} \}.
\]
Now 
\[ \alpha^\Up = \big ( \bigcup { A \mid (A,B) \in \mathcal{H}} \} \big )^{\Up \DOWN \Up} 
= \big ( \bigcup { A \mid (A,B) \in \mathcal{H}} \} \big )^{\Up} = \bigcup \{ A^\Up \mid (A,B) \in \mathcal{H}\} .
\]
For each $B \in \wp(U)^\UP$, $B = B^{\Down \UP}$. We have
\[ \beta^\Down = 
\big (\bigcup \{B \mid (A,B) \in \mathcal{H} \} \big )^\Down = 
\big (\bigcup \{B^{\Down \UP} \mid (A,B) \in \mathcal{H} \} \big )^\Down 
= \big (\bigcup \{B^{\Down} \mid (A,B) \in \mathcal{H} \} \big )^{\UP\Down}. \]
Since $B^\Down \in \varphi[\mathcal{T}]$ for each $(A,B) \in \mathcal{T}$
and $\varphi[\mathcal{T}]$ is closed with respect to arbitrary unions, we have that
$\bigcup \{ B^\Down \mid (A,B) \in \mathcal{H} \}$ belongs to $\varphi[\mathcal{T}]$. This means that 
\[ \beta^\Down = 
\big (\bigcup \{B^{\Down} \mid (A,B) \in \mathcal{H} \} \big )^{\UP\Down} =
\bigcup \{B^{\Down} \mid (A,B) \in \mathcal{H} \} =
\bigcup \{A^{\Up} \mid (A,B) \in \mathcal{H} \} = \alpha^\Up. 
\]
Thus, $\bigvee \mathcal{H} \in \mathcal{T}$ and hence, $\mathcal{T}$
is a complete sublattice of $\mathrm{DM(RS)}$.
\end{proof}

Our next result can be viewed as a generalization of Corollary~\ref{cor:BoolComplement}. 

\begin{corollary}\label{cor:Boolean}
Let $R$ be a reflexive relation on $U$ and $\mathcal{T} \subseteq \mathcal{C}$ a nonempty set closed with respect to ${\sim}$. If $\mathcal{T}$ is a  sublattice of $\mathrm{DM(RS)}$, 
then it is a Boolean lattice.
\end{corollary}

\begin{proof}
Recall from Lemma~\ref{Lemma:A-set} that $\mathcal{A}$ is closed
with respect to set-theoretical complement and $(\mathcal{C},{\leq} ,{\sim} )$ and $(\mathcal{A},{\subseteq}, {^c})$ 
are isomorphic pseudo-Kleene posets via the map $\varphi$, according to Proposition \ref{Prop:isomorphism}. If $\mathcal{T}$ is a sublattice
of $\mathrm{DM(RS)}$, then $\varphi[T]$ is a sublattice of $\wp(U)$. In addition, if
$\mathcal{T}$ is closed with respect to $\sim$, then $\varphi[\mathcal{T}]$ is closed 
under $c$. Thus, $\varphi[T]$ is a Boolean lattice, and so is $\mathcal{T}$.
\end{proof}

\section{Central and exact elements}
\label{Sec:CentralElements}

In this section, we consider the central elements of $\mathrm{DM(RS)}$
and their relationship with the exact rough sets. The central elements of a bounded lattice $L$
have  an important role, because they correspond exactly to the 
direct decompositions of $L$ \cite[Theorem III.4.1]{Gratzer1998}.

In the well-studied case of rough sets induced by equivalences, the central elements are the exact rough sets. 
This observation can be found in 
\cite{GeWa92}, where a precise structure theorem was
given by showing that $\mathrm{RS}$ is always isomorphic to a product
of chains of two and three elements.

Here, we prove that the set of exact rough sets of $\mathrm{DM(RS)}$ equals
the intersection of the central elements of $\mathrm{DM}(\mathrm{RS})(R)$ and $\mathrm{DM}(\mathrm{RS})(R^{-1})$;
see below for the description of the notation. This means that if $R$ is a tolerance, then central elements 
coincide with exact rough sets. We end this section by noting that
if $R$ is a tolerance induced by an irredundant covering or a quasiorder,
then central elements are the same as sharp, complemented, or 
exact rough sets.

We begin this section by recalling some definitions and facts from \cite{Gratzer1998}.
Let $L$ be a lattice. An element $a \in L$ is called \emph{neutral} if 
\[ (a \wedge x) \vee (x \wedge y) \vee (y \wedge a) = 
 (a \vee x) \wedge (x \vee y) \wedge (y \vee a)
\]
for all $x,y \in L$. It is known that $a \in L$ is neutral if
and only if the sublattice generated by $a$, $x$, and $y$ is distributive for any $x,y \in L$.

Let $L$ be a bounded lattice. An element is said to be \emph{central} if 
it is complemented and neutral.
The set of central elements of $L$ is denoted by $\mathrm{Cen}(L)$.
As noted in \cite{Gratzer1998}, $0,1 \in \mathrm{Cen}(L)$
and $\mathrm{Cen}(L)$ forms a Boolean sublattice of $L$.
Thus, the complement of a central element is both unique and also a central element.

In \cite[Theorem 4.13]{Maeda1970}, the following conditions
are shown to be equivalent:
\begin{enumerate}[({M}1)]
\item $a \in  \mathrm{Cen}(L)$;
\item there exists an element $a'$ such that
\[ x = (x \wedge a) \vee (x \wedge a') = (x \vee a) \wedge (x \vee a')\] 
for every $x \in L$.
\end{enumerate}

We can now write a characterization of central elements in terms of the
operation ${\sim}$.

\begin{lemma}\label{Lemma:centinv}
Let $(L,\vee ,\wedge ,\sim ,0,1)$ be a bounded  lattice with an antitone
involution. If $a\in L$ is sharp, then the following are equivalent:
\begin{enumerate}[\rm (i)]
\item $a \in \mathrm{Cen}(L)$;
\item $x = (x \wedge a) \vee (x \wedge {\sim} a)$ for all $x\in L$;
\item $x = (x \vee a) \wedge (x \vee {\sim }a)$ for all $x\in L$.
\end{enumerate}
\end{lemma}

\begin{proof}
Let $a$ be sharp. We know that $a$ is complemented and ${\sim} a$ is a complement of $a$. 
Suppose that $a \in \mathrm{Cen}(L)$ and $x \in L$. 
Because the sublattice generated by $a$, ${\sim} a$, and $x$ is distributive, 
we have
\begin{align*}
x &= x \wedge 1 = x \wedge (a \vee {\sim} a)
= (x \wedge a) \vee (x \wedge {\sim} a) \\
\intertext{and} 
x &= x \vee 0 = x \vee (a \wedge {\sim} a) = 
(x \vee a) \wedge (x\vee  {\sim} a).
\end{align*}
Thus, (ii) and (iii) hold. Because the lattice $L$ is self-dual, 
(ii) and (iii) are equivalent. Indeed, assume that (ii) is true.
Then, for $y \in L$,  
\begin{align*}
{\sim} y & = {\sim}((y \wedge a) \vee (y \wedge {\sim} a)) =
 {\sim} (y \wedge a) \wedge {\sim} (y \wedge {\sim} a) \\
  & = ({\sim} y \vee {\sim} a) \wedge ({\sim}y \vee {\sim} {\sim} a) =
   ({\sim}y \vee a) \wedge ({\sim} y \vee {\sim} a). 
\end{align*}
Let $x \in L$. Set $y := {\sim} x$. Then, $x = {\sim} y$ and
(iii) holds for $x$. Thus, (ii) implies (iii). Similarly, we can show
that (iii) implies (ii). 

If we combine (ii) and (iii), we have that
there exists an element ${\sim}a$ such that 
\[ x = (x \wedge a) \vee (x \wedge {\sim} a)
= (x \vee a) \wedge (x\vee  {\sim} a).
\]
By the equivalence of (M1) and (M2), we obtain that $a \in \mathrm{Cen}(L)$.
\end{proof}

Let $R$ be a reflexive relation on $U$. Suppose that 
$(A,B) \in \mathrm{DM(RS)}$ is central. Since $(A,B)$ is (uniquely) complemented, Proposition~\ref{Prop:UniqueComplement} 
says ${\sim}(A,B) = (B^c,A^c)$ is
its unique complement which is also a central element. Moreover, $A^\Up = B^\Down$ and
$(A,B) \in \mathrm{RS}$ by  Lemma~\ref{lem:simple_sharp}.

Our next result characterizes central elements among sharp elements.

\begin{proposition} \label{Prop:CentRS}
Let $R$ be a reflexive relation on $U$ and let 
$(A,B)$ be sharp in $\mathrm{DM(RS)}$.
Then $(A,B)$ is central if and only if 
$X = ( (X\cap A) \cup (X \cap B^c) )^{\Up \DOWN}$
for any $X\in \wp(U)^\DOWN$.
\end{proposition}

\begin{proof}
Because ${\sim} (A,B)=(B^{c},A^{c})$, in view of Lemma~\ref{Lemma:centinv}, $(A,B)$ is a central element if and only if
\begin{equation} \tag{$\star$} \label{eq:central}
(X,Y) = ((X,Y)\wedge (A,B))\vee ((X,Y) \wedge (B^{c},A^{c}))
\end{equation}
for each $(X,Y)\in \mathrm{DM(RS)}$. This equation is equivalent to
\[
(X,Y) = ((X\cap A,(Y\cap B)^{\Down \UP })\vee 
((X\cap B^{c}),(Y\cap A^{c})^{\Down \UP }).
\]
We can expand the lattice operations further and obtain
\begin{align*}
(X,Y) &= \big ( (X\cap A)\vee (X\cap B^{c}),
(Y\cap B)^{\Down \UP } \cup (Y\cap A^{c})^{\Down \UP } \big)  \\
& = \big ( ((X\cap A)\cup (X\cap B^{c}))^{\Up \DOWN},
(Y\cap B)^{\Down \UP} \cup (Y\cap A^{c})^{\Down \UP} \big ).
\end{align*}
Thus, we have that $(A,B)$ is a central if and only if
\begin{align}
X & = ( (X\cap A)\cup (X\cap B^{c}) )^{\Up \DOWN } 
\label{equation:first} \\
\intertext{and}
Y &= (Y \cap B)^{\Down \UP} \cup (Y \cap A^{c})^{\Down \UP}
\label{equation:second}
\end{align}
for each $(X,Y) \in \mathrm{DM(RS)}$.

We prove that \eqref{equation:second} is satisfied for every 
$Y \in \wp(U)^{\UP}$. Because $(A,B)$ is sharp, $B^\Down = A^\Up$. 
We have that
\begin{align*}
(Y\cap B)^{\Down \UP } \cup (Y\cap A^{c})^{\Down \UP } 
& = ((Y^{\Down }\cap B^{\Down })\cup (Y^{\Down }\cap A^{c\Down}))^{\UP} \\
& = (Y^{\Down }\cap (B^{\Down }\cap A^{c\Down }))^\UP \\
& = (Y^{\Down }\cap (A^{\Up }\cap A^{\Up c}))^{\UP } \\
& =Y^{\Down \UP} = Y.
\end{align*}

It is now clear that if $(A,B)$ is central, then \eqref{equation:first} 
holds for any $X \in \wp(U)^\DOWN$. On the other hand, suppose
that $(X,Y) \in \mathrm{DM(RS)}$. Then $X \in \wp(U)^\DOWN$ and
$Y \in \wp(U)^\UP$. Assume that \eqref{equation:first} is true for $X$.
Because $Y$ satisfies \eqref{equation:second} trivially, by the above 
considerations, \eqref{eq:central} is satisfied for $(X,Y)$. Thus,
$(A,B)$ is central, and the proof is completed.
\end{proof}

Let $R$ be a reflexive relation on $U$. A rough set 
$(A,B) \in \mathrm{RS}$ is called \emph{exact} whenever $A=B$. Observe that in this case $A^{\UP} = A^{\Up } = A = A^{\DOWN}=A^{\Down}$. 
Indeed, $(A,B) \in \mathrm{RS} \subseteq \mathrm{DM(RS)}$ implies 
$A^{\Up \UP}\subseteq B$. As $R$ is reflexive, we obtain
$A\subseteq A^{\Up} \subseteq A^{\Up \UP} \subseteq B=A$ and 
$A \subseteq A^{\UP} \subseteq A^{\Up \UP} \subseteq B=A$, whence we
get $A=A^{\Up}=A^{\UP}$. 

As it is noted in \cite{Jarvinen2007}, $({^\UP},{^\Down})$ and $({^\Up},{^\DOWN})$ form order-preserving Galois 
connections on $(\wp(U),\subseteq)$.
Therefore, $A = A^\UP$ and $A = A^\Up$ give
$A \subseteq A^{\UP \Down} = A^\Down \subseteq A$ and 
$A \subseteq A^{\Up \DOWN} = A^\DOWN \subseteq A$, that is,
$A = A^\DOWN = A^\Down$. Because $B=A$, we have $B^{\UP} = B^{\Up} = B = B^{\DOWN} = B^{\Down}$. 
As these relations yield $B^{\Down} = B = A =A^{\Up}$, in view of Lemma \ref{Prop:SharpCharacterization}, 
an exact rough set $(A,B)$ is a
complemented and sharp element of $\mathrm{DM(RS)}$.

In addition, since $A^{\Down \UP} \subseteq A \subseteq A^{\UP \Down}$ and
$A^{\DOWN \Up} \subseteq A \subseteq A^{\Up \DOWN}$ and $R$ is reflexive, 
$A = A^\UP$ implies $A \subseteq A^\Down \subseteq A$, that is, $A = A^\Down$. Similarly, 
$A^\Down = A$ yields $A^\UP \subseteq A \subseteq A^\UP$. We can prove a similar relationship between $A^\Up$ and $A^\DOWN$.
Thus, the following equivalences hold for all $A \subseteq U$:
\begin{align}
A^\UP = A &\iff A^\Down = A; \label{eq:fix1} \\
A^\Up = A &\iff A^\DOWN = A. \label{eq:fix2}
\end{align}

\begin{lemma}\label{Lem:exact_equalities}
If $R$ is a reflexive relation, the following are equivalent:
\begin{enumerate}[\rm (i)]
\item $(A,A)$ is a rough set;
\item $A^{\DOWN} = A^{\UP}$;
\item $A^{\Down} = A^{\Up}$;
\item $A^{\UP} = A^{\Up} = A$;
\item $A^{\DOWN} = A^{\Down} = A$.
\end{enumerate}
\end{lemma}

\begin{proof}
Suppose that (i) holds. Then $(A,A)$ is an exact rough set, and
(ii)--(v) follow immediately from the previous observations.

Since $R$ is reflexive, (ii) is equivalent to $A^{\DOWN}= A = A^{\UP}$.
Thus, $(A,A)$ is a rough set and (i) holds. 
If (ii) holds, then $A^\DOWN = A = A^\UP$, and  by \eqref{eq:fix1}
and \eqref{eq:fix2}, we have $A^\Down = A = A^\Up$. Hence, (ii) implies (iii). Similarly, we can see that (iii) implies (ii).

Suppose that (iv) holds. Then $A^\UP = A$ and $A^\Up = A$
imply $A^\DOWN = A$ and $A^\Down = A$ by \eqref{eq:fix1}
and \eqref{eq:fix2}. Thus, (v) holds. Analogously, we can prove
that (v) implies (iv).

Clearly, (ii) and (iii) together imply (iv) and (v).
Conversely, (iv) implies (ii) and from (v) we infer (iii)
by using \eqref{eq:fix1} and \eqref{eq:fix2}.
\end{proof}

\begin{lemma}\label{Lemma:exact}
If $R$ is a reflexive relation, then any exact rough set is a central
element of $\mathrm{DM(RS)}$.
\end{lemma}

\begin{proof}
Let $(A,A)$ be an exact rough set. Then $A^{\Down}=A^{\Up }$ holds by Lemma~\ref{Lem:exact_equalities}(iii). This means that $(A,A)$ is
sharp. Let $X \in \wp(U)^{\DOWN}$. Then $X^{\Up \DOWN} =X$ and 
\[ X^{\Up \DOWN} = (X\cap (A\cup A^{c})^{\Up \DOWN} =
((X \cap A)\cup (X \cap A^{c}))^{\Up \DOWN}. \] 
Thus, we obtain $X = ((X\cap A) \cup (X\cap A^{c}))^{\Up \DOWN }$, 
that is, \eqref{equation:first} holds. Hence,
in view of Proposition \ref{Prop:CentRS}, 
$(A,A)$ is a central element of $\mathrm{DM(RS)}$.
\end{proof}

An \emph{equivalence} is a reflexive, symmetric and transitive binary relation. 
For an equivalence $E$ on $U$, we denote by $x/E$ the \emph{equivalence class}
of $x \in U$, that is, the set of all elements $y$ related to $x$. The family of
all $E$-classes is denoted by $U/E$.

For any binary relation $R$ on $U$, its transitive closure 
$R^+$ is defined as the least transitive relation on $U$ containing 
$R$. It is well-known that
\[  R^{+} = \bigcup _{i \geq 1} R^{i}, \]
where $R^i$ is the $i$-th power of $R$, defined inductively by
$R^{1}=R$ and for $i > 1$, $R^{i+1} = R\circ R^{i}$, 
where $\circ$ denotes the composition of relations.
Let $R^e$ denote the smallest equivalence containing $R$. Obviously,
$R^e = (R\cup R^{-1})^+$.  

The following result is presented in \cite[Prop.~6.1.8]{Jarv99} for tolerances, 
here we extend it to reflexive relations.

\begin{proposition}\label{Prop:equivalence_class}
Let $R$ be a reflexive relation on $U$. Then $(A,A)$ is exact 
if and only if $A$ is a union of $R^e$-classes.
\end{proposition}

\begin{proof}
Assume that $A= \bigcup \mathcal{H}$ for some $\mathcal{H} \subseteq U/{R^e}$. Every $X \in \mathcal{H}$ is an $R^e$-class, that is, 
$X = a/R^e$ for some $a \in X$.
If $x \in X^\UP$, then there is $y \in X$ such that $(x,y) \in R \subseteq R^e$. Since $(y,a) \in R^e$, 
we get $(x,a) \in R^a$ by the transitivity of $R^e$. Thus,  $x \in X$ and $X^{\UP} \subseteq X$. The inclusion $X \subseteq X^\UP$ holds by the reflexivity 
of $R$. Similarly, we can show that $X^{\Up} = X$ for all 
$X \in \mathcal{H}$. Now,
\[ A^{\UP} = \bigcup_{X \in \mathcal{H}} X^\UP
= \bigcup_{X \in \mathcal{H}} X = A \text{ \ and \ }
A^{\Up} = \bigcup_{X \in \mathcal{H}} X^\Up = 
\bigcup_{X \in \mathcal{H}} X = A. 
\]
In view of Lemma~\ref{Lem:exact_equalities}, $(A,A)$ is exact.

Conversely, let $(A,A)$ be an exact rough set. We prove
that $A = \bigcup \mathcal{H}$ for $\mathcal{H} = \{ x / R^e \mid  x \in A\}$.
Because $a \in a / R^e$, $A \subseteq \bigcup \mathcal{H}$.
Suppose $a \in \bigcup \mathcal{H}$. Then $a \in b/R^e$ for some $b \in A$. Because $(a,b) \in R^e = (R \cup R^{-1})^+$,
there exists a sequence $z_{0},z_{1},...,z_{n}\in U$ such that $a=z_{0}$, $b = z_{n}$ and $(z_{i-1},z_{i}) \in R\cup R^{-1}$ 
for each $1\leq i\leq n$.

As $(A,A)$ is exact, we have $A^\UP = A$ and $A^\Up = A$. This
implies that each $z_{n},z_{n-1},\ldots,z_{1}, z_0$ belongs to
$A$. In particular, $a \in A$. This proves $\bigcup \mathcal{H}
\subseteq A$.
\end{proof}

Let $L$ be a lattice with a least element $0$. 
The lattice $L$ is \emph{atomistic}, if any element of $L$ is the join of
atoms below it. It is well known (see e.g. \cite{Gratzer1998}) that a 
complete Boolean lattice is atomistic if and only if it is completely distributive,
that is, a complete lattice in which arbitrary joins distribute over arbitrary meets. 
It is known that any complete sublattice of $\wp(U)$ is atomistic.

Let $E$ be an equivalence on $U$. We say that a set $X \subseteq U$
is \emph{saturated} by $E$ if $X = \bigcup \mathcal{H}$ for some 
$\mathcal{H} \subseteq U/E$.
Let us denote by $\mathrm{Sat}(E)$ the set of all $E$-saturated sets.
In the rough set literature, sets saturated by $E$ are also called \emph{$E$-definable}.
Recall that a set $X \subseteq U$ is \emph{$E$-definable} whenever $X^\DOWN = X^\UP$.
It is known \cite{Steinby97} that saturated sets form a complete atomistic Boolean lattice such that 
\[ \bigvee \mathcal{H} = \bigcup \mathcal{H} \quad \text{and} \quad
\bigwedge \mathcal{H} = \bigcap \mathcal{H} \]
for $\mathcal{H} \subseteq \mathrm{Sat}(E)$. The complement of 
$X \in \mathrm{Sat}(E)$ is $X^c$ and $U/E$ is the set of atoms of $\mathrm{Sat}(E)$.

It is clear that the set of $E$-exact rough sets is 
$\{(X,X) \mid X \in \mathrm{Sat}(E)\}$ and the map $X \mapsto (X,X)$ is a trivial order-isomorphism between $\mathrm{Sat}(E)$ 
and $E$-exact rough sets. This means
that $E$-exact sets form a complete atomistic Boolean lattice with respect to the 
coordinatewise inclusion such that 
\[ \bigvee \{ (A,A) \mid A \in \mathcal{H}\} =
\big (\bigcup \mathcal{H}, \bigcup \mathcal {H} \big )  
\quad \text{and} \quad
\bigwedge \{ (A,A) \mid A \in \mathcal{H}\} =
\big ( \bigcap \mathcal{H}, \bigcap \mathcal{H} \big )  
\]
for all $\mathcal{H} \subseteq \mathrm{Sat}(E)$.
The complement of an exact set $(A,A)$ is $(A^c,A^c)$
and the set of atoms of exact sets is $\{ (a/E,a/E) \mid a \in U\}$.

The following corollary is clear by Proposition~\ref{Prop:equivalence_class} and
the fact $R^e = (R^{-1})^e$.

\begin{corollary}\label{cor:R_and_inverse}
Let $R$ be reflexive. The exact rough sets induced by $R$, $R^{-1}$, and $R^e$ are the same.
\end{corollary}

Corollary~\ref{cor:R_and_inverse} means that $R$-exact sets form a complete 
atomistic Boolean lattice  such that 
\[ \bigvee \{ (A,A) \mid A \in \mathcal{H}\} =
\big (\bigcup \mathcal{H}, \bigcup \mathcal {H} \big )  
\quad \text{and} \quad
\bigwedge \{ (A,A) \mid A \in \mathcal{H}\} =
\big (\bigcap \mathcal{H}, \bigcap \mathcal{H} \big )  
\]
for all $\mathcal{H} \subseteq \mathrm{Sat}(R^e)$.
The set of atoms is $\{ (a/R^e,a/R^e) \mid a \in U\}$. The
complement of $(A,A)$ is ${\sim}(A,A) = (A^c,A^c)$
for $A \in \mathrm{Sat}(R^e)$. It is easy to observe that
$R$-exact rough sets form a complete sublattice of
$\mathrm{DM(RS)}$. For instance, if $\mathcal{H} \subseteq \mathrm{Sat}(R^e)$,
then the join in $\mathrm{DM(RS)}$ is
\[ \bigvee \{ (A,A) \mid A \in \mathcal{H} \} =
\big ( \big (\bigcup \mathcal{H}\big)^{\Up \DOWN}, \bigcup \mathcal{H} \big)
= \big (\bigcup \mathcal{H}, \bigcup \mathcal{H} \big).
\]
The latter equality follows from the fact that $\bigcup \mathcal{H}$
is $R^e$-saturated.

A reflexive and transitive binary relation is called a \emph{quasiorder}.
We know \cite{JRV09} that if $R$ is a quasiorder on $U$, then $\mathrm{RS}$ is a complete sublattice of 
$\wp(U) \times \wp(U)$. Thus, $\mathrm{DM(RS)} = \mathrm{RS}$ is bounded and (completely) distributive. 

\begin{proposition}\label{Prop:cent1}
If $R$ is a quasiorder on $U$, then for any $(A,B) \in \mathrm{RS}$ the following are equivalent:
\begin{enumerate}[\rm (i)]
    \item $(A,B)$ is sharp;
    \item $(A,B)$ is complemented.
    \item $(A,B)$ is a central element of the lattice $\mathrm{RS}$.
    \item $(A,B)$ is an exact rough set.
\end{enumerate}
\end{proposition}

\begin{proof}
We have that (i) and (ii) are equivalent according to 
Lemma~\ref{Prop:SharpCharacterization}. 
Since $\mathrm{RS}$ is bounded and distributive, (ii) implies (iii). Trivially, (iii) implies (ii).

(iii)$\Rightarrow $(iv): Assume that $(A,B)=(Z^\DOWN,Z^\UP) \in RS$ is 
a central element. Then $(A,B)$ is complemented and sharp, and $A^{\Up} = B^{\Down}$. This yields $Z^{\DOWN \Up} = Z^{\UP \Down}$. As $R$
is a quasiorder, we have $Z^{\DOWN \Up}=Z^\DOWN$ and $Z^{\UP \Down}
=Z^\UP$. Hence, $A=Z^\DOWN = Z^{\DOWN \Up} = Z^{\UP \Down} = Z^{\UP}=B$. This means that $(A,B)$ is an exact rough set. 

Finally, (iv) implies (iii) according to Lemma~\ref{Lemma:exact}.
\end{proof}

\begin{lemma}\label{lem:cent2} Let $R$ be a reflexive
relation on $U$. If $(A,B)$ is a central element of $\mathrm{DM(RS)}$, 
then $A = B^\DOWN$.
\end{lemma}

\begin{proof}
As $(A,B) \in \mathrm{DM(RS)}$, $A^{\Up \UP} \subseteq B$. This gives
$A \subseteq A^\Up \subseteq  (A^\Up)^{\UP \DOWN}=(A^{\Up \UP})^\DOWN \subseteq B^\DOWN$.
Since $(A,B)$ is central,
\[  X = ((X\cap A) \cup (X\cap B^{c}))^{\Up \DOWN}\]
for each $X\in \wp(U)^\DOWN$ by Proposition~\ref{Prop:CentRS}. 
Set $X:=B^\DOWN$. Then 
$X \cap A = B^\DOWN \cap A = A$ and $X \cap B^{c} = 
B^{\DOWN} \cap B^{c} \subseteq B \cap B^c = \emptyset$. Hence, 
we obtain $B^\DOWN = (A \cup \emptyset)^{\Up \DOWN} = A^{\Up \DOWN} = A$.
The last equality follows from the fact that $A \in \wp(U)^\DOWN$.
\end{proof}

Let $\mathcal{L} := \mathrm{DM}(\mathrm{RS})(R)$ and 
$\mathcal{L}^{-1} := \mathrm{DM}(\mathrm{RS})(R^{-1})$, that is,
$\mathcal{L}$ and $\mathcal{L}^{-1}$ are the smallest completions
of the ordered sets of rough sets determined by $R$ and $R^{-1}$,
respectively.

\begin{proposition}\label{Prop:cent2}
If $R$ is a reflexive relation on $U$, then the set of exact rough 
sets of $R$ equals 
$\mathrm{Cen}(\mathcal{L}) \cap \mathrm{Cen}(\mathcal{L}^{-1})$.
\end{proposition}

\begin{proof}
Let $(A,A)$ be an exact rough set defined by $R$. Then, in view of
Corollary \ref{cor:R_and_inverse}, $(A,A)$ is also an exact rough set of 
$R^{-1}$. By Lemma \ref{Lemma:exact}, $(A,A)$ is at the same time a central element of $\mathcal{L}$ and of $\mathcal{L}^{-1}$.

Conversely, assume $(A,B) \in \mathrm{Cen}(\mathcal{L}) \cap \mathrm{Cen}(\mathcal{L}^{-1})$. 
Then $(A,B)$ is sharp with respect to $R$ and
$R^{-1}$, and $(A,B) \in \mathrm{RS}$. By Lemma~\ref{lem:cent2},
$A^\Up = B^\Down = A$ and $A^\UP = B^\DOWN = A$. We get
$B = B^{\Down \UP} = A^\UP = A$, proving that $(A,B)$ is exact.
\end{proof}

A \emph{tolerance relation} is a reflexive and symmetric binary relation.

\begin{corollary}\label{Cor:excent}
If $R$ is a tolerance relation on $U$, then the central elements of $\mathrm{DM(RS)}$ coincide to the exact rough sets induced
by $R$.
\end{corollary}

\begin{proof}
Since $R=R^{-1}$, the lattices $\mathcal{L}$ and $\mathcal{L}^{-1}$
are the same. Hence, the central elements of $\mathrm{DM(RS)}$ coincide  with the exact rough sets of $R$.
\end{proof}

A collection $\mathcal{H}\subseteq \wp(U)$ of nonempty
subsets of $U$ is called a \emph{covering} of $U$ if 
$\bigcup \mathcal{H}=U$. $\mathcal{H}$ is an \emph{irredundant covering}, if
$\mathcal{H}\backslash \{X\}$ is not a covering of $U$ for any $X\in
\mathcal{H}.$ We proved in  \cite{JarRad2019} that a tolerance relation 
$R$ is induced by an irredundant covering if and only if for any $(a, b) \in R $, there exists an element $c \in U$ such that $a,b \in R(c)$ and
$R(c)$ is a block. A \emph{block} is a subset of $U$ whose elements are 
$R$-related to each other. We  also proved that in this case $\mathrm{RS}$ 
is a complete, completely distributive lattice. Because in a distributive complete
lattice the complemented elements and the central elements of the lattice are the same,
we obtain the following corollary.

\begin{corollary}\label{Cor:cent}
Let $R$ be a tolerance induced by an irredundant covering of $U$. Then for any $(A,B) \in $
$\mathrm{RS}$ the following are equivalent:
\begin{enumerate}[\rm (i)]
    \item $(A,B)$ is sharp;
    \item $(A,B)$ is complemented;
    \item $(A,B)$ is a central element of the lattice $\mathrm{RS}$;
    \item $(A,B)$ is an exact rough set.
\end{enumerate}
\end{corollary}

In our next example, we show that for reflexive relations, there are cases in which complemented---and thus sharp, elements are not exact.

\begin{example} \label{exa:non_exact_complemented}
Let us consider the pseudo-Kleene algebra of Remark~\ref{rem:Chajda}. Its Hasse diagram is in Figure~\ref{fig:twoFigs}(a). 
Let us find the central elements. By definition, the central elements are complemented.

The complemented elements are $(\emptyset,\emptyset)$, $(\{a\},\{a,b\})$, $(\{c\},\{b,c\})$ and $(U,U)$. 
The elements $(\emptyset,\emptyset)$ and $(U,U)$ are exact. They are central, because exact elements are always central by 
Proposition~\ref{Prop:cent2}.  

The element $(\{a\},\{a,b\})$ is not central, because the sublattice generated by 
$(\{a\},\{a,b\})$,  $(\{c\},U)$ and $(\{b,c\},U)$ is 
\[ \{(\emptyset, \{a,b\}),  (\{a\},\{a,b\}), (\{c\},U), (\{b,c\},U), (U,U)\}, \]
which is not distributive. 

The element $(\{c\},\{b,c\})$ cannot be central, because its complement $(\{a\},\{a,b\})$ is not central. 
Recall that the complement of a central element is central.
\end{example}

\section{Basic concepts of Brouwer--Zadeh lattices}
\label{sec:BZ-lattices}

In this section, we recall from the literature the essential facts 
about Brouwer--Zadeh lattices. They are pseudo-Kleene 
algebras provided with an additional negation $\Neg$ resembling an intuitionistic
negation. The two negations are connected by the identity 
${\sim} \Neg x = \Neg \Neg x$. In terms of the two negations, one can
define a necessity and a possibility operators. There is a close connection between the 
set $\mathcal{N}$ of clopen elements of these operators and the operation $\Neg$.

A \emph{Brouwer--Zadeh lattice}  is an algebra $(L,\vee ,\wedge ,{\sim},\Neg ,0,1)$
such that $(L,\vee ,\wedge ,\sim ,0,1)$ is a pseudo-Kleene algebra and for all $a,b \in L$ the following conditions hold:
\begin{enumerate}[\rm ({BZ}1)]
\item $a \wedge \Neg a=0$;
\item $a \leq \Neg \Neg a$;
\item $a\leq b$  implies $\Neg b \leq \Neg a$;
\item ${\sim} \Neg a = \Neg \Neg a$.
\end{enumerate}
Brouwer--Zadeh lattices are also called simply \emph{BZ-lattices}.
A paraorthomodular BZ-lattice is called a \emph{PBZ}-lattice.

Let $(L,\vee ,\wedge ,{\sim},\Neg ,0,1)$ be a BZ-lattice.
In \cite{GLP17} it is noted that the following facts hold for all $a \in L$:
\begin{enumerate}[\rm ({BZ}1)]
\setcounter{enumi}{4}
\item $\Neg a \leq  {\sim} a$;
\item $\Neg \Neg \Neg a = \Neg a$;
\item $\Neg a$ and ${\sim} \Neg a$ are complements of each other.
\end{enumerate}

As in \cite{Cattaneo84}, we define a pair of mappings on $L$ by setting
\[ \Diamond a := \Neg \Neg a \quad \text{and} \quad \Box a := \Neg {\sim} a.\]

\begin{lemma} \label{lem:closure_interior}
Let $(L,\vee ,\wedge ,{\sim},\Neg ,0,1)$ be a BZ-lattice and $a,b \in L$.
\begin{enumerate}[\rm (i)]
\item $\Box a \leq a \leq \Diamond a$;
\item $a \leq b$ implies $\Box a \leq \Box b$ and $\Diamond a \leq \Diamond b$;
\item $\Box \Box a = \Box a$ and $\Diamond \Diamond a = \Diamond a$;
\item $\Box \Diamond a = \Diamond a$ and $\Diamond \Box a = \Box x$;
\item ${\sim} \Diamond a = \Box {\sim} a$ and ${\sim} \Box a = \Diamond {\sim} a$. 
\end{enumerate}
\end{lemma}

\begin{proof}
(i) By (BZ2), $a \leq \Neg\Neg a = \Diamond a$. By (BZ5),
$\Box a = \Neg {\sim} a \leq {\sim}{\sim} a = a$.  

(ii) Let $a \leq b$. By (BZ3), $\Neg b \leq \Neg a$ and $\Diamond a = \Neg\Neg a \leq \Neg\Neg b = \Diamond b$.
We have also  ${\sim} b \leq {\sim} a$, which gives $\Box a = \Neg {\sim} a \leq \Neg {\sim} b = \Box b$.

(iii) By (BZ4), $\Box \Box a = \Neg {\sim} \Neg {\sim} a = \Neg \Neg \Neg {\sim} a = \Neg {\sim} a = \Box a$. We get
$\Diamond \Diamond a = \Neg \Neg \Neg \Neg a = \Neg \Neg a = \Diamond a$ by (BZ6).

(iv) We have by (BZ4) and (BZ6) that $\Box \Diamond a = \Neg {\sim} \Neg \Neg a = \Neg\Neg\Neg\Neg a = \Neg\Neg a = \Diamond a$.
Similarly, $\Diamond \Box a = \Neg \Neg \Neg {\sim} a = \Neg {\sim} a = \Box a$ by (BZ6).

(v) Fact (BZ4) gives ${\sim \Diamond} a = {\sim}\Neg\Neg a = \Neg \Neg \Neg a = \Neg a = \Neg {\sim} {\sim} a = \Box {\sim a}$.
Finally, ${\sim} \Box a = {\sim} \Neg {\sim} a = \Neg \Neg {\sim} a = \Diamond {\sim} a$.
\end{proof}
Note that $\Box$ and $\Diamond$ can be defined in terms of each other, that is, for all $a \in L$,
\[ \Box a = {\sim} \Diamond {\sim} a \quad \text{and} \quad \Diamond a = {\sim} \Box {\sim} a.\]

Let $P$ be an ordered set. Then a function $c \colon P \to P$ is called a {\em closure operator}
on $P$, if for all $a,b \in P$, $a \leq c(a)$,  $a \leq b$ implies $c(a) \leq c(b)$, and
$c(c(a)) = c(a)$. An element $a \in P$ is called {\em closed} if $c(a) = a$.
The set of all $c$-closed elements of $P$ is denoted by $P_c$.

Let $c$ be a closure operator on an ordered set $P$. Then the following holds:
\begin{enumerate}[\rm ({C}1)]
\item $P_c = \{c(a) \mid a \in P\}$;
\item if $k$ is a closure operator on $P$ with $P_c = P_k$, then $c = k$;
\item $c(x) = \bigwedge_P \{a \in P_c \mid x \leq a\}$ for any $x \in P$;
\item if $\bigvee S$ exists in $P$, then
$\bigvee S$ exists in $P_c$ and $\bigvee_{P_c} S = c(\bigvee_P S)$
for all $S \subseteq P_c$;
\item if $\bigwedge S$ exists in $P$, then $\bigwedge S$ exists
in $P_c$ and $\bigwedge_{P_c} S = \bigwedge_P S$
for all $S \subseteq P_c$.
\end{enumerate}

An \emph{interior operator} is defined dually. More precisely, a function $i \colon P \to P$ is
an interior operator on $P$, if for all $a,b \in P$, $i(a) \leq a$,  $a \leq b$ implies $i(a) \leq i(b)$, and
$i(i(a)) = i(a)$. An element $a \in P$ is called {\em open} if $i(a) = a$. 

\begin{remark} \label{rem:open_closed_N}
Let $(L,\vee ,\wedge ,{\sim},\Neg, 0,1)$ be a BZ-lattice. By cases (i)--(iii) of Lemma~\ref{lem:closure_interior} $\Diamond$ is a closure
operator and $\Box$ is an interior operator. Lemma~\ref{lem:closure_interior}(iv)
means that the set of $\Diamond$-closed elements and $\Box$-open are equal and is denoted by $\mathcal{N}$, that is,
\[ \mathcal{N} = \{ \Box a \mid a \in L\} = \{ \Diamond a \mid a \in L\}. \]
Note that elements of $\mathcal{N}$ can be viewed as \emph{clopen}---they are
both $\Diamond$-closed and $\Box$-open at the same time.
Obviously, $\mathcal{N}$ is a sublattice of $L$ and if $L$ is a complete lattice,
then $\mathcal{N}$ is its complete sublattice.

It is also easy to see that
\[ \mathcal{N} = \{ \Neg x \mid x \in L\}. \]  
Indeed,  if $x \in \mathcal{N}$, then $x = \Diamond a = \Neg \Neg a$ for some $a \in L$.
Conversely, consider $\Neg a$ for some $a \in L$. Now $\Diamond \Neg a = \Neg \Neg \Neg a = \Neg a$.
Thus, $\Neg a \in \mathcal{N}$.
 
In \cite{Giuntini2016}, $\Diamond$-closed elements are called $\Diamond$-sharp.
As mentioned in \cite{Giuntini2016}, each element of $\mathcal{N}$ is sharp 
in the sense we considered in Section~\ref{Sec:basic}. Indeed, let $a \in \mathcal{N}$. Then $a = \Neg b$ for some $b \in L$.
According to (BZ7), $a$ has a complement ${\sim} \Neg b = {\sim} a$. Note that this means that $\mathcal{N}$ is closed under $\sim$.
It is now clear that $(\mathcal{N}, \vee, \wedge, {\sim}, \Neg, 0, 1)$ is a subalgebra of  $(L,\vee ,\wedge ,{\sim},\Neg ,0,1)$.

An element $a$ is called \emph{Brouwer-sharp} if $a \vee \Neg a = 1$.
Each $\Diamond$-closed elements is Brouwer-sharp, because $a = \Diamond a$
implies $a \vee \Neg a = \Neg \Neg a \vee \Neg a = {\sim}\Neg a \vee \Neg\Neg\Neg a  
= {\sim}\Neg a \vee {\sim} \Neg\Neg a 
= {\sim}(\Neg a \wedge \Neg \Neg a) = {\sim} 0 = 1$.

Additionally, each Brouwer-sharp element is sharp. Indeed,
$1 = a \vee \Neg a \leq a \vee {\sim} a$. In \cite{Giuntini2016} it is
proved that in PBZ-lattices, $\Diamond$-closed elements
and Brouwer-sharp elements coincide. However, there are PBZ-lattices in which all $\Diamond$-closed elements are not sharp.

Note also that in \cite{Giuntini2016} 'sharp elements' are called 'Kleene sharp elements' in distinction with 
$\Diamond$-sharp and Brouwer-sharp elements.
\end{remark}

An \emph{ortholattice} is an algebra $(L,\vee,\wedge,{^\bot},0,1)$ such that
$(L,\vee,\wedge,0,1)$ is a bounded lattice and ${^\bot}$ is an 
antitone involution satisfying
$x \vee x^\bot = 1$ and $x \wedge x^\bot = 0$.
Let $(L,\vee ,\wedge ,{\sim},0,1)$ be a pseudo-Kleene algebra such that the underlying lattice $L$ is complete. 
If $K$ is a complete pseudo-Kleene subalgebra of $L$ such that
$(L,\vee ,\wedge ,{\sim},0,1)$ forms an ortholattice, then we say
that $K$ is a \emph{complete subortholattice} of $L$. Note that
in such a case, the orthocomplementation of $K$ is always assumed to be
the Kleene negation $\sim$ of $L$.

\begin{proposition} \label{Prop:OrtholatticeFromBZ}
Let $(L,\vee ,\wedge ,{\sim},\Neg ,0,1)$ be a BZ-lattice such that the underlying lattice $L$ is complete. 
Then $\mathcal{N}$ is a complete subortholattice of  
the pseudo-Kleene algebra $(L,\vee ,\wedge ,{\sim},0,1)$.
\end{proposition}

\begin{proof}
By Remark~\ref{rem:open_closed_N}, $(\mathcal{N},\vee,\wedge,{\sim},0,1)$ is a
subalgebra of $(L,\vee,\wedge, {\sim},0,1)$ as a pseudo-Kleene algebra.
Further, each $a \in \mathcal{N}$ has a complement ${\sim}a$. Thus,
$(\mathcal{N}, \vee, \wedge, {\sim}, 0, 1)$ is an ortholattice.

If $L$ is a complete lattice, then, by (C4) and (C5), $\bigvee_L S$
and $\bigwedge_L S$ belong to $\mathcal{N}$ for all $S \subseteq \mathcal{N}$.
Thus, $\mathcal{N}$ is a complete sublattice of $L$.
\end{proof}

Let $L$ be a complete lattice. If $\mathcal{S}$ is a complete sublattice of $L$, then it is known that the
map $\Diamond^\mathcal{S}$ on $L$ defined by
\[ \Diamond^{\mathcal{S}} x = \bigwedge \{ a \in \mathcal{S} \mid x \leq a\} \]
is a closure operator such that the set of $\Diamond^{\mathcal{S}}$-closed elements is $\mathcal{S}$.
Similarly, $\Box^{\mathcal{S}}$ defined by 
\[ \Box^{\mathcal{S}} x = \bigvee \{ a \in \mathcal{S} \mid a \leq x\} \]
is an interior operator on $L$ such that the set of open elements is $\mathcal{S}$.

We can now write an 'opposite statement' of Proposition~\ref{Prop:OrtholatticeFromBZ}.

\begin{proposition} \label{Prop:OppositeOrtholatticeP}
Let $(L,\vee ,\wedge ,{\sim},0,1)$ be a pseudo-Kleene algebra such that the underlying lattice $L$ is
complete. If $\mathcal{N}$ is a complete subortholattice of $L$, then 
$(L,\vee ,\wedge ,{\sim},\Neg,0,1)$ is a BZ-lattice in which the
operation $\Neg$ is defined for all $x \in L$ by
\[ \Neg x :=  \Box^{\mathcal{N}} {\sim} x = {\sim} \Diamond^{\mathcal{N}} x.\] 
\end{proposition}

\begin{proof}
Let us denote here $\Box^{\mathcal{N}}$ and $\Diamond^{\mathcal{N}}$ simply by $\Box$ and $\Diamond$, respectively.
We show that $\Neg$ satisfies conditions  (BZ1)--(BZ4).

(BZ1) Because $\Diamond x$ is in $\mathcal{N}$, we have $\Diamond x \wedge {\sim} \Diamond x = 0$. 
Now ${\sim} \Neg x = \Neg \Neg x = \Diamond x$ implies 
${\sim}\Diamond x = {\sim}{\sim} \Neg x = \Neg x$. We obtain
$x \wedge \Neg x \leq \Diamond x \wedge {\sim} \Diamond x = 0$.

(BZ2) $\Neg \Neg x = \Neg {\sim} \Diamond x = \Box {\sim}{\sim} \Diamond x = \Box \Diamond x = \Diamond x \geq x$.

(BZ3) $x \leq y$ implies ${\sim} y \leq {\sim} x$ and $\Neg y = \Box {\sim} y 
\leq \Box {\sim} x = \Neg x$.

(BZ4) By definition, $\Neg \Neg x = \Diamond x$. Because $\Neg x = {\sim} \Diamond x,
{\sim} \Neg x = {\sim} {\sim} \Diamond x = \Diamond x$. Therefore, 
$\Neg \Neg x = {\sim}\Neg x$.
\end{proof}

A \emph{BZ*-lattice} is a BZ-lattice $(L, \vee, \wedge, {\sim} ,\Neg ,0,1)$ that satisfies for all $a \in L$ the condition
\begin{enumerate}[(BZ8)]
\item[(BZ8)] $\Neg (a \wedge  {\sim} a) \leq \Neg a \vee \Neg {\sim a}$.
\label{eq:condstar}
\end{enumerate}
A \emph{PBZ*-lattice} is a paraorthomodular BZ*-lattice.

A BZ-lattice whose underlying lattice distributive is called
a \emph{distributive BZ-lattice}. Any distributive BZ*-lattice is a PBZ*-lattice. Indeed, assume that
$a \leq b$ and ${\sim} a \wedge b = 0$. Then $a \wedge b = a$ and $a$ and $b$
are sharp. Using distributivity, we have (cf. \cite[Lemma~2.3]{Giuntini2016}):
\[
    a = (a \wedge b) \vee ({\sim} a \wedge b) =
    b \wedge ((a \wedge b) \vee {\sim} a) = 
    b \wedge (a \vee {\sim}a) = b \vee 1 = b.
\]

Let us recall from \cite{Gratzer1998} some facts related to 
pseudocomplemented lattices and Stone algebras.
A \emph{pseudocomplemented lattice} is a bounded lattice $L$
equipped with a unary operation $^*$ characterized by the property:
\[
a \wedge x=0 \iff  x\leq a^*.   
\]
The element $a^*$ is called the pseudocomplement of $a$. 
By definition, for any $a, b \in L$, $a \leq b$ implies $ b^* \leq a^*$.  
Also $a \leq a^{**}$, because $a^* \wedge a^{**} = 0$. 

A distributive pseudocomplemented lattice satisfying the \emph{Stone identity}
\[
x^* \vee x^{**} = 1
\]
is called a \emph{Stone algebra}. In any Stone algebra $L$, the equality
\[
(a\wedge b)^* = a^* \vee b^*
\]
is satisfied for all $a,b\in L$. In addition, the elements 
$\{ a^* \mid a \in L \}$ form a Boolean sublattice of $L$ in which
$^*$ is the complement operation.

A distributive pseudo-Kleene algebra is a \emph{Kleene algebra}.
A \emph{pseudocomplemented Kleene algebra} $(L,\vee ,\wedge,{\sim},^{\ast},0,1)$ is
a Kleene algebra $(L,\vee ,\wedge,{\sim},0,1)$ with a pseudocomplementation.
Furthermore, a pseudocomplemented Kleene algebra 
satisfying the Stone identity  is called a \emph{Kleene--Stone algebra}.

\begin{proposition}\label{prop:BZstar)}
A Kleene--Stone algebra such that its complemented and sharp elements coincide is a PBZ*-lattice.
\end{proposition}

\begin{proof}
Trivially, the reduct  $(L,\vee ,\wedge,{\sim}, 0,1)$ of a Kleene--Stone
algebra $(L,\vee ,\wedge,{\sim},^{\ast},0,1)$ is a pseudo-Kleene algebra. 
We show that $^*$ satisfies (BZ1)--(BZ4). Because
$a \wedge a^* = 0$, (BZ1) holds. Since $a\leq a^{**}$, condition
(BZ2) is satisfied. As $a\leq b$ implies $b^* \leq a^*$, (BZ3) also holds.

Because $L$ is a Stone lattice, $a^{\ast}$ has a unique 
complement $a^{\ast \ast}$ for each $a \in L$. As by our
assumption $a^*$ is a sharp element, ${\sim} a^*$ is
also a complement of $a^*$. Thus, we obtain 
${\sim} a^* = a^{**}$. Hence, condition (BZ4) holds and we have proved that
$(L,\vee ,\wedge,{\sim},^{\ast},0,1)$ is a BZ-lattice.

Since $L$ is a Stone lattice, 
\[ (a \wedge {\sim} a)^* = a^* \vee ({\sim} a)^*\] 
Thus, (BZ8) holds, and we have a BZ*-lattice. 

Because each Stone algebra is distributive by definition, 
this BZ*-lattice is a PZB*-lattice.
\end{proof}

\section{PBZ and PBZ*-lattices on the completion of rough sets}
\label{sec:PBZ*-lattices}

Our main result of this section shows that there is a bijective
correspondence between atomistic complete Boolean 
sublattices of $\mathrm{DM(RS)}$ and  PBZ-lattices defined on $\mathrm{DM(RS)}$.
Because $\{(\emptyset,\emptyset),(U,U)\}$ forms a trivial such a 
sublattice, we can always define at least one PBZ-lattice on $\mathrm{DM(RS)}$.
We will also consider PBZ-lattices defined by quasiorders and tolerances induced by irredundant coverings.

Let $E$ be an equivalence on $U$. Then, $\mathrm{RS}$ is a complete sublattice of $\wp(U) \times \wp(U)$.
It is proved by \cite{PomPom88} that $\mathrm{RS}$ is a Stone algebra in which $(A,B)^* = (B^c,B^c)$.
Therefore, $\mathrm{RS}$ forms a Kleene--Stone algebra. Note that $\sim$ is defined as usual, that is,
${\sim}(A,B) = (B^c,A^c)$ for all $(A,B) \in \mathrm{RS}$. By Proposition~\ref{Prop:cent1}, the exact, sharp,
complemented and central elements of $\mathrm{RS}$ coincide.
Applying Proposition~\ref{prop:BZstar)}, we can write the following.

\begin{proposition} \label{prop:PBZ_equivalence}
For an equivalence $E$ on $U$, the rough sets defined by $E$ form a PBZ*-lattice
$(\mathrm{RS},\vee,\wedge,{\sim},{^*},(\emptyset,\emptyset), (U,U)$.
\end{proposition}

Let $R$ be a reflexive binary relation on $U$. As we have already shown in Section~\ref{Sec:basic},
$\mathrm{DM(RS)}$ forms a paraorthomodular complete pseudo-Kleene algebra. 
In addition, the sharp elements coincide with complemented elements $\mathcal{C}$.
We say that an equivalence $E$ on $U$ \emph{extends} $R$ (or $E$ is an equivalence
\emph{extending} $R$) if $R \subseteq E$. Next, we show how to define PBZ-lattice on
$\mathrm{DM(RS)}$ in terms of any equivalence relation $E$ extending $R$. 
Because we have here two relations, we denote by $X^\uparrow$ and $X^\downarrow$
the upper and lower approximations of a set $X$ determined by $E$.
The approximation determined by $R$ are denoted as usual, that is, by $X^\DOWN$ and $X^\UP$.

We will need the following lemma.

\begin{lemma} \label{lem:E_extending_R}
Let $E$ be an equivalence extending a reflexive relation $R$ on $U$.
Then for all $X \subseteq U$,
\[
X^{\uparrow \UP} = X^{\uparrow \DOWN} = X^\uparrow \quad \text{and} \quad
X^{\downarrow \UP} = X^{\downarrow \DOWN} = X^\downarrow.\]
\end{lemma}

\begin{proof}
Because $R$ is reflexive, $X^{\downarrow \DOWN} \subseteq X^\downarrow$. Let $x \in X^\downarrow$ and $y \in R(x)$.
Suppose for contradiction that $y \notin X^\downarrow$. Then there is $z \notin X$ such that $(y,z) \in E$.
Because $(x,y) \in R \subseteq E$, we have $(x,z) \in E$. The fact $x \in X^\downarrow$ gives $z \in X$,
a contradiction. Hence, $y \in X^\downarrow$ and $R(x) \subseteq X^\downarrow$. Thus, $x \in X^{\downarrow \DOWN}$ and 
$X^{\downarrow \DOWN} = X^\downarrow$.

We have $X^\downarrow \subseteq X^{\downarrow \UP}$. Because $R \subseteq E$, we get $X^{\downarrow \UP} \subseteq X^{\downarrow \uparrow} = X^\downarrow$.
Similarly, $X^{\uparrow\DOWN} \subseteq X^\uparrow$. To prove the converse, assume $x \notin X^{\uparrow \DOWN}$. Then,
there is $y \in R(x)$ such that $y \notin X^\uparrow$, that is, $y/E \cap X = \emptyset$. Now
$(x,y) \in R \subseteq E$ gives $x/E = y/E$. Thus, $x/E \cap X = y/E \cap X = \emptyset$ and $x \notin X^\uparrow$.
\end{proof}

\begin{proposition} \label{prop:ExtendingEquivalence}
Let $R$ be a reflexive relation on $U$ and let $E$ be an equivalence extending $R$.
We obtain a PBZ-lattice $(\mathrm{DM(RS)},\vee ,\wedge ,\sim ,\Neg,(\emptyset ,\emptyset ),(U,U))$ by setting
\begin{equation*}
\Neg (A,B) := (B^{c \downarrow},  B^{c \downarrow})  
\end{equation*}
for any $(A,B)\in \mathrm{DM(RS)}$.
\end{proposition}

\begin{proof}
We know that $\mathrm{DM(RS)}$ forms a paraorthomodular pseudo-Kleene algebra.
We check first that the map $\Neg$ is well-defined.  Let $(A,B) \in \mathrm{DM(RS)}$. Then, by Lemma~\ref{lem:E_extending_R},
\[ \Neg (A,B) = ( B^{c \downarrow}, B^{c \downarrow}) =  ( B^{c \downarrow \DOWN}, B^{c \downarrow \UP}) .\]
Thus, $\Neg (A,B)$ belongs to $\mathrm{RS} \subseteq \mathrm{DM(RS)}$.

Next we verify the properties of the operation $\Neg$. Let $(A,B), (C,D) \in \mathrm{DM(RS)}$.

\smallskip\noindent%
(BZ1) Because $A\subseteq B$, we obtain $B^{c}\subseteq A^{c}$ and $B^{c \downarrow}  \subseteq B^c \subseteq A^c$. 
We have $A \cap B^{c \downarrow} \subseteq A \cap A^c = \emptyset$ and 
$(B \cap B^{c \downarrow})^{\Down \UP} \subseteq B \cap B^{c \downarrow} \subseteq B \cap B^c = \emptyset$.
Note that the map $X \mapsto X^{\Down \UP}$ is an interior operator on $\wp(U)$; see \cite{Jarvinen2007}. We have that
\[ (A,B) \wedge \Neg (A,B) = (A,B) \wedge ( B^{c \downarrow}, B^{c \downarrow} )  = 
( A \cap B^{c \downarrow}, ( B \cap  B^{c \downarrow} )^{\triangledown \blacktriangle })
\subseteq (\emptyset, \emptyset).
\]
\smallskip\noindent%
(BZ2) If $(A,B)\leq (C,D)$, then $B \subseteq D$ implies $D^{c \downarrow} \subseteq B^{c \downarrow}$, whence we obtain
\[ \Neg (C,D) = (  D^{c \downarrow},  D^{c \downarrow}) \leq ( B^{c \downarrow}, B^{c \downarrow}) = \Neg (A,B). \]

\smallskip\noindent%
(BZ3) By definition,
\[ \Neg \Neg (A,B) = \Neg ( B^{c \downarrow} , B^{c \downarrow} ) =    ( B^{c \downarrow c \downarrow} , B^{c \downarrow c \downarrow} )  =
 ( B^{\uparrow c c \downarrow} , B^{\uparrow c c \downarrow} ) =  ( B^{\uparrow  \downarrow} , B^{\uparrow  \downarrow} ) 
=  ( B^{\uparrow} , B^{\uparrow} ) \geq (A,B), \]
because $B^\uparrow \supseteq B \supseteq A$.

\smallskip\noindent%
(BZ4) By direct computation, ${\sim} \Neg (A,B) = {\sim} ( B^{c \downarrow}, B^{c \downarrow}) =( B^{c \downarrow c}, B^{c \downarrow c}) = (B^\uparrow, B^\uparrow)$.
Hence, ${\sim} \Neg (A,B)= \Neg\Neg (A,B)$.
\end{proof}

\begin{example} \label{exa:quasiorders}
Let $R$ be a quasiorder on $U = \{a,b,c\}$ such that $R(a) = \{a,b\}$, $R(b) = \{b\}$, $R(c) = \{c\}$. Let us denote in this example
sets as sequences of their elements, like $\{a,b\}$ is denoted by $ab$. The Hasse diagram of the lattice
$\mathrm{RS}$ is given in Figure~\ref{fig:quasiorders}.

Note that the complemented elements are $(ab,ab)$ and $(c,c)$. These two are also
the sharp, exact, and central elements of $\mathrm{RS}$.
\begin{figure}
    \centering
    \includegraphics[width=50mm]{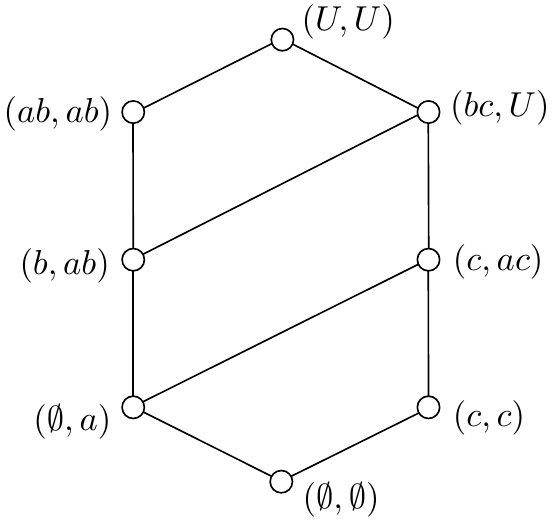}
    \caption{The lattice $\mathrm{RS}$ of Example~\ref{exa:quasiorders}\label{fig:quasiorders}}
\end{figure}
There are two equivalences extending $R$: (i) $R^e$, which has the equivalence classes $\{a,b\}$ and $\{c\}$, and (ii) the universal relation $U \times U$.
Let us denote the corresponding negations by $\Neg_1$ and $\Neg_2$, respectively. Their values are given in the following table:
{\small
\[\begin{array}{r|cccccccc}
       (A,B) & (\emptyset,\emptyset) & (\emptyset,a) & (b,ab) &(ab,ab) & (c,c) & (c,ac) 
       & (bc,U) & (U,U)  \\ \hline
\Neg_1 (A,B) & (U,U) & (c,c)  & (c,c) & (c,c)  &(ab,ab) & (\emptyset,\emptyset)& (\emptyset, \emptyset) & (\emptyset,\emptyset)  \\
\Neg_2 (A,B) & (U,U) &(\emptyset,\emptyset)&(\emptyset,\emptyset)&(\emptyset,\emptyset)&(\emptyset,\emptyset)&(\emptyset,\emptyset)&(\emptyset,\emptyset)&(\emptyset,\emptyset)
\end{array}
\]
}
\end{example}

As can be seen in Example~\ref{exa:quasiorders}, the sets $\{ \Neg_i (A,B) \mid (A,B) \in \mathrm{DM(RS)}\}$, for $1\leq i \leq 2$, 
are (complete) subortholattices of $\mathrm{RS}$. Actually, we may present the following correspondences.

\begin{proposition} \label{prop:correspondence}
Let $R$ be a reflexive relation. If 
$(\mathrm{DM(RS)}, \vee ,\wedge ,\sim ,\Neg ,(\emptyset,\emptyset),(U,U))$ is a PBZ-lattice,
then $\mathcal{N}^\Neg := \{\Neg (A,B) \mid (A,B) \in \mathrm{DM(RS)}\}$ is a  complete subortholattice of $\mathrm{DM(RS)}$.
Conversely, given a complete subortholattice $\mathcal{N}$ of $\mathrm{DM(RS)}$, the formula
\[ \Neg^\mathcal{N} (A,B) := \bigvee \{ (X,Y) \in \mathcal{N} \mid X \cap B = \emptyset \text{ and } Y \cap A = \emptyset \}\]
defines a PBZ-lattice $(\mathrm{DM(RS)}, \vee ,\wedge ,\sim ,\Neg^\mathcal{N} ,(\emptyset,\emptyset),(U,U))$. 
\end{proposition}

\begin{proof}
By Proposition~\ref{Prop:OrtholatticeFromBZ}, $\mathcal{N}^\Neg$ is a complete subortholattice of $\mathrm{DM(RS)}$. 
On the other hand, if $\mathcal{N}$ is a complete subortholattice of $\mathrm{DM(RS)}$,
using Proposition~\ref{Prop:OppositeOrtholatticeP}, we can define the operation $\Neg^\mathcal{N}$ by the formula
\begin{align*}
  \Neg^\mathcal{N} (A,B) &=  \Box^{\mathcal{N}} {\sim} (A,B) \\
  & = \bigvee \{ (X,Y) \in \mathcal{N} \mid (X,Y) \leq {\sim}(A,B) \} \\
  & = \bigvee \{ (X,Y) \in \mathcal{N} \mid (X,Y) \leq (B^c,A^c) \} \\
      & = \bigvee \{ (X,Y) \in \mathcal{N} \mid X \subseteq B^c \text{ and } Y \subseteq A^c \} \\
  &= \bigvee \{ (X,Y) \in \mathcal{N} \mid X \cap B = \emptyset \text{ and } Y \cap A = \emptyset \}. \qedhere
\end{align*} 
\end{proof}

\begin{corollary}\label{cor:N_is_Boolean}
Let $R$ be a reflexive relation. If $(\mathrm{DM(RS)}, \vee ,\wedge ,\sim ,\Neg ,(\emptyset,\emptyset),(U,U))$ is a PBZ-lattice,
then $\mathcal{N}^\Neg = \{\Neg x \mid x \in \mathrm{DM(RS)}\}$ forms an atomistic complete Boolean sublattice of $\mathrm{DM(RS)}$ 
in which $\sim$ is the complement operation.
\end{corollary}

\begin{proof}
By Remark~\ref{rem:open_closed_N}, $\mathcal{N}^\Neg$ is a complete sublattice of $\mathrm{DM(RS)}$. 
Because $\mathcal{N}^\Neg$ is closed under $\sim$, it is complemented, and thus a Boolean lattice by Corollary~\ref{cor:Boolean}.
Because $\mathcal{N}^\Neg$ is order-isomorphic to a completely
distributive lattice $\varphi[N^\Neg]$, $\mathcal{N}^\Neg$ is completely distributive and hence atomistic. 
\end{proof}

By Proposition~\ref{prop:correspondence} and Corollary~\ref{cor:N_is_Boolean},
we can write the following theorem, because the atomistic complete Boolean 
sublattice induced by the operator $\Neg^\mathcal{N}$ is $\mathcal{N}$ itself, and, analogously, the operator induced by $\mathcal{N}^\Neg$ is $\Neg$.

\begin{theorem}\label{thm:correspondence}
Let $R$ be a reflexive relation. The correspondence between atomistic complete Boolean 
sublattices of $\mathrm{DM(RS)}$ and  PBZ-lattices on $\mathrm{DM(RS)}$ is bijective.
\end{theorem}

Our next theorem states that PBZ-lattices determined by quasiorders or tolerances induced by an irreduncant covering are such that the
operation $\Neg$ is always determined by an equivalence extending these relations.

\begin{theorem}\label{Thm:PBZ_characterization}
If $R$ is a quasiorder on $U$ or a tolerance induced by an irredundant covering of $U$,
then $(\mathrm{DM(RS)}, \vee, \wedge, {\sim},\Neg, (\emptyset,\emptyset),(U,U))$ is a
PBZ-lattice if and only if there exists an equivalence $E$ extending $R$ and $\Neg (A,B)$ equals $(B^{c \downarrow}, B^{c \downarrow})$
for all $(A,B) \in \mathrm{RS}$.
\end{theorem}

\begin{proof}
Assume $R$ is a quasiorder on $U$ or a tolerance induced by an irredundant covering of $U$. In both cases, $R$ is reflexive.
Let $E$ be an equivalence extending $R$. Then by Proposition~\ref{prop:ExtendingEquivalence}, we may define a PBZ-lattice on
$\mathrm{RS}$ by setting $\Neg (A,B) = ( B^{c \downarrow},  B^{c \downarrow})$ for any $(A,B)\in \mathrm{RS}$.

On the other hand, we know by Corollary~\ref{cor:N_is_Boolean} that $\mathcal{N} = \{\Neg x \mid x \in \mathrm{RS}\}$ forms an
atomistic complete Boolean sublattice of $\mathrm{RS}$. Because each element of $\mathcal{N}$ is complemented, $\mathcal{N}$
consists of exact elements of $\mathrm{RS}$. Let $\mathrm{At}$ denote the set of atoms of $\mathcal{N}$. Then
$\mathrm{At} = \{(A,A) \mid A \in \mathcal{H}\}$ for some $\mathcal{H} \subseteq \mathrm{Sat}(R^e)$ according to
Proposition~\ref{Prop:equivalence_class}. As the meet of two different atoms of $\mathcal{N}$ is $(\emptyset,\emptyset)$,
we have $A_1 \cap A_2 \neq \emptyset$ for all $A_1,A_2 \in \mathcal{H}$. Since the greatest element of $\mathcal{N}$ is
$(U,U)$, we have $(U,U) = \bigvee \mathrm{At} = (\bigcup \mathcal{H}, \bigcup \mathcal{H})$. This means that $\mathcal{H}$
forms a partition of $U$.

Let $E$ be the equivalence corresponding to $\mathcal{H}$. We have that any $E$-class is a union of $R^e$-classes.
This means that $(x,y) \in R^e$ implies $(x,y) \in E$ and $R \subseteq R^e \subseteq E$. Thus, $E$ extends $R$.

For $(C,D) \in \mathrm{RS}$,
\begin{align*}
 \Neg (C,D) &= \Box^\mathcal{N} {\sim}(C,D) = \bigvee \{(X,X) \in \mathcal{N} \mid (X,X) \leq {\sim}(C,D) \} \\
 &= \bigvee \{ (X,X) \in \mathcal{N} \mid (X,X) \leq (D^c,C^c) \} \\
 &= \bigvee \{ (A,A) \in \mathrm{At} \mid (A,A) \leq (D^c,C^c) \} \\
 &= \bigvee \{ (A,A)  \mid A \in U/E \text{ and } (A,A) \leq (D^c,C^c) \}.
\end{align*}
As $C \subseteq D$ gives $D^c \subseteq C^c$, we have $(A,A) \leq (D^c,C^c)$ if and only if $A \leq D^c$.
We obtain
\begin{align*}
\Neg (C,D) &= \bigvee \{(A,A) \mid A \in U/E \text{ and } A \subseteq B^c \} \\
&= \big (\{\bigcup \{ A \in U/E \mid A \subseteq B^c \},\{\bigcup \{ A \in U/E \mid A \subseteq B^c \}\big )
= (B^{c\downarrow},B^{c\downarrow} ).  \qedhere
\end{align*}
\end{proof}

\begin{example}
Let us consider the pseudo-Kleene algebra of Figure~\ref{fig:twoFigs}(a). Then,
\[\mathcal{C} = \{ (\emptyset,\emptyset), (\{a\},\{a,b\}), (\{c\},\{b,c\}), (U,U)\}. \]
For $(A,B) \in \mathrm{RS}$, we have
\[ \Diamond^\mathcal{C}(A,B) = \bigwedge \{ (X,Y) \in \mathcal{C} \mid (A,B) \leq (X,Y) \} .\]
Because $\Neg (A,B) = {\sim} \Diamond^\mathcal{C}(A,B)$, we can write the following table:
{\small
\[
\begin{array}{lll}
  (A,B) & \Diamond^\mathcal{C} (A,B) & \Neg(A,B) \\ \hline
  (\emptyset,\emptyset) & (\emptyset,\emptyset) & (U,U) \\
  (\emptyset,\{a\})     & (\{a\},\{a,b\})  & (\{c\}, \{b,c\}) \\
  (\{c\}, \{b,c\})      & (\{c\}, \{b,c\}) & (\{a\},\{a,b\}) \\
  (\emptyset,\{a,b\})   & (\{a\},\{a,b\})  & (\{c\}, \{b,c\}) \\
  (\{c\}, U)            & (U, U) & (\emptyset,\emptyset) \\
  (\{a\},\{a,b\})       & (\{a\},\{a,b\})  & (\{c\}, \{b,c\}) \\
  (\{b,c\}, U)            & (U, U) & (\emptyset,\emptyset) \\
  (U, U)            & (U, U) & (\emptyset,\emptyset)          
\end{array}
\]
}

Because $R$ is not a quasiorder or a tolerance induced by an irredundant covering, the elements $\Neg(A,B)$ are not
necessarily exact sets.
\end{example}

\begin{lemma} \label{lem:bz8}
Let $E$ be an equivalence extending a quasiorder $R$ on $U$. The operation $\Neg$ on $\mathrm{RS}$
defined by $\Neg (A,B) = (B^{c \downarrow},B^{c \downarrow})$ satisfies \eqref{eq:condstar} if and only if
$(A \cup B^c)^\downarrow \subseteq A^\downarrow \cup B^{c \downarrow}$ for all $(A,B) \in \mathrm{RS}$.
\end{lemma}

\begin{proof} Condition \eqref{eq:condstar} means that
\[ \Neg ((A,B) \wedge (B^c,A^c)) \leq \Neg(A,B) \vee \Neg (B^c,A^c).\]
Because $R$ is a quasiorder, $\mathrm{RS}$ is a complete sublattice of $\wp(U) \times \wp(U)$. We have
\[ \Neg ((A,B) \wedge (B^c,A^c)) = \Neg(A \cap B^c, B \cap A^c) = ((B\cap A^c)^{c \downarrow},(B\cap A^c)^{c \downarrow})
= ((A \cup B^c)^{\downarrow}, (A \cup B^c)^{\downarrow})\]
and
\[  \Neg(A,B) \vee \Neg (B^c,A^c) = (B^{c \downarrow},B^{c \downarrow}) \vee (A^\downarrow,A^\downarrow) =
(A^\downarrow \cup B^{c \downarrow},A^\downarrow \cup B^{c \downarrow}).\]
Condition \eqref{eq:condstar} is thus equivalent to $(A \cup B^c)^{\downarrow} \subseteq A^\downarrow \cup B^{c \downarrow}$.
\end{proof}

For our following theorem, we need to recall a couple of our older results. Let $R$ be a quasiorder. It is proved in
\cite[Proposition 4.2]{JPR13} that a pair $(A,B) \in \wp(U)^\DOWN \times \wp(U)^\UP$ belongs to $\mathrm{RS}$ if and only
if $A \subseteq B$ and $\mathcal{S} \subseteq A \cup B^c$.
Recall that $\mathcal{S}$ denotes the set of all singletons.
In addition, in \cite[Theorem 6.4]{JRV09} it is showed that $\mathrm{R}$S is a Stone lattice if and only if $R^{-1} \circ R = R^e$.

\begin{theorem} \label{thm:PBZ_star}
Let $R$ be a quasiorder on $U$. Then $(\mathrm{RS},\vee ,\wedge ,{\sim} ,\Neg,(\emptyset ,\emptyset ),(U,U))$
is a distributive PBZ*-lattice if and only if the operation $\Neg$ is defined 
for any $(A,B) \in \mathrm{DM(RS)}$ by
\begin{equation}
\Neg (A,B) := ( B^{c \downarrow}, B^{c \downarrow})
\end{equation}
in terms of $R^e$. If $R\circ R^{-1} = R^e$, then $\Neg$ is the pseudocomplementation.
\end{theorem}
\begin{proof}
As we already mentioned, for any quasiorder, $\mathrm{RS}$ is always a distributive lattice.
Let $E$ be an equivalence extending $R$. In view of Theorem~\ref{Thm:PBZ_characterization}, we can define a PBZ-lattice
on $\mathrm{RS}$ if and only if $\Neg (A,B) = (B^{c\downarrow},B^{c\downarrow})$. 

Let us now consider the equivalence $R^e$. For clarity, we denote in this proof 
by $X^\Downarrow$ the lower approximation of $X \subseteq U$ determined by $R^e$. 
We prove that $(A\cup B^{c})^\Downarrow \subseteq A^\Downarrow \cup B^{c\Downarrow}$.
Let us assume that this does not hold. Thus, there exists an element $x \in (A\cup B^{c})^\Downarrow$ such that
$x \notin (A^\Downarrow \cup B^{c\Downarrow})$. Therefore, $x/R^e \nsubseteq A$ and $x/R^e \nsubseteq B^c$.
Because $A \subseteq B$, we have $A \cap B^c = \emptyset$. Now $x/R^e \nsubseteq A$ implies $x/R^e \cap B^c \neq \emptyset$, because
$x / R^e \subseteq A \cup B^c$. Similarly, $x / R^e \nsubseteq B^c$ implies $x / R^e \cap A \neq \emptyset$.
Thus, there exist elements $a \in A$ and $b \in B^c$ such that $(a,b) \in R^e$.

Between any elements $a \in A$ and $b \in B^c$ there can be several \emph{paths} $z_0, z_1,\ldots,z_n$ such that $z_0 = a$, $z_n = b$ and
$(z_{i-1},z_{i}) \in R \cup R^{-1}$ for $1 \leq i \leq n$, where $n$ means the \emph{length} of the path. Because the
length of each path is a nonnegative integer, there exists at least one path with the minimal length $n$ connecting
$a$ and $b$. We say that the \emph{distance} between the sets $A$ and $B^c$ is the minimal length between all possible paths between
some $a \in A$ and $b \in B^c$.

Let the distance between $A$ and $B^c$ be $n$. In addition, let $z_0 \in A$, $z_n \in B^c$ and $z_0, z_1,\ldots,z_n$ be the
corresponding connecting path. Because $A \cap B^c = \emptyset$, $n = 0$ is impossible. Similarly,
$n=1$ would imply that $(z_{0},z_{1}) \in R\cup R^{-1}$. Thus, $(u,v) \in R$ or $(v,u)\in R$.
As $(A,B) \in \mathrm{RS}$, we have $A=X^\DOWN$ and $B=X^\UP$ for some $X\subseteq U$.
Now, $z_0 \in A=X^\DOWN$ and $(z_0,z_1)\in R$ imply $z_1 \in R(z_0)\subseteq X \subseteq X^\UP = B$.
We get $z_1 \in B^{c}\cap B=\emptyset $, a
contradiction. Similarly, $(z_1,z_0)\in R$ and $z_0 \in A$ yield $z_1 \in A^\UP \subseteq X^\UP = B$. Again, 
$z_1 \in B^{c}\cap B=\emptyset$, a contradiction.

Thus, $n \geq 2$. As the path $z_{0},z_{1},\dots,z_{n}$ connecting $A$ and $B^c$ is of the minimal length, we must have 
$z_{1} \notin A$ and $z_{1} \notin B^{c}$.
Indeed, having $z_1 \in A$ would imply that $z_1,z_2,\ldots, z_n$ is a path connecting
$A$ and $B$ shorter than the path $z_0, z_1 \dots, z_n$  of the minimum length. Similarly, $z_1 \in B^c$ would imply that $z_0, z_1$ is a
path of length $1$ connecting $A$ and $B^c$, which case we just managed to exclude. 

Because $z_1 \notin A \cup B^c$, we have $x/R^e \nsubseteq A \cup B^{c}$, contradicting the assumption $x \in (A\cup B^{c})^\Downarrow$.
This proves that for any  $x \in (A\cup B^{c})^\Downarrow$, we must have  $x\in A^\Downarrow \cup B^{c\Downarrow}$, that is,
$(A\cup B^{c})^\Downarrow \subseteq A^\Downarrow \cup B^{c\Downarrow}$.

\medskip%

Let $E$ be an equivalence extending $R$. Note that $R \subseteq E$ is equivalent to
$R^e \subseteq E$. We assume $R^e \subset E$ and show that there is $(A,B) \in \mathrm{RS}$
such that $(A \cup B^c)^\downarrow \nsubseteq A^\downarrow \cup B^{c \downarrow}$, where
$X^\downarrow$ denotes the lower approximation of $X \subseteq U$ determined by $E$.

Because $R^e \subset E$, there exists an $E$-class $H$ that is a union of at least
two $R^e$-classes. Let $K$ be the union of $R^e$-classes of the elements 
in $\mathcal{S} \setminus H$, that is,
$K = \bigcup \{ z/R^e \mid z \in \mathcal{S} \setminus H \}$.
Notice that by definition  
$K \subseteq \bigcup \{ z/E \mid z \in \mathcal{S} \setminus H \} \subseteq H^c$.

Let $x \in H$. We set
\[ A := x / R^e \cup K \quad \text{and} \quad B := x / R^e \cup H^c.\]
It is clear that $A$ and $B$ belong to $\mathrm{Sat}(R^e)$. Thus, $A^\Downarrow = A$ and
$B^\Downarrow = B$. By Lemma~\ref{lem:E_extending_R}, $A^\DOWN = A^{\Downarrow \DOWN} = A^\Downarrow = A$.
Thus, $A \in \wp(U)^\DOWN$. Similarly, we can show $B \in \wp(U)^\UP$.
By definition, $A \subseteq B$.

We prove that $A \cap \mathcal{S} = B \cap \mathcal{S}$. Clearly, $A \cap \mathcal{S}
\subseteq B \cap \mathcal{S}$. Suppose that $a \in B \cap \mathcal{S}$. If $a \in x/R^e$, then
obviously $a \in A \cap \mathcal{S}$. If $a \in H^c$, then $a \in \mathcal{S} \setminus H$
and $a \in a/R^e$ imply $a \in K$ and $a \in A$. Hence, $a \in A \cap \mathcal{S}$ and
$A \cap \mathcal{S} = B \cap \mathcal{S}$ holds.  Thus,
$\mathcal{S} = (B \cap \mathcal{S}) \cup (B^c \cap \mathcal{S})  =
 (A \cap \mathcal{S}) \cup (B^c \cap \mathcal{S})
\subseteq A \cup B^c$ and $(A,B) \in \mathrm{RS}$.

Now,
\begin{align*}
  A \cup B^c &= (x/R^e \cup K) \cup (x/R^e \cup H^c)^c = (x/R^e \cup K) \cup ((x/R^e)^c \cap H) \\
  &= (x/R^e \cup K \cup (x/R^e)^c) \cap (x/R^e \cup K \cup H) = U  \cap (x/R^e \cup K \cup H) \\
  &= x/R^e \cup K \cup H = K \cup H.
\end{align*}
The last equality follows from the fact that $x/R^e \subset x/E = H$. We have that $x/E 
= H \subseteq K \cup H = A \cup B^c$ and $x \in (A \cup B^c)^\downarrow$. On the other 
hand, $K \subseteq H^c$ means $K \cap H = \emptyset$. Now $x / R^e \subseteq x/E$ imply
$x/E = H \nsubseteq x/R^e \cup K = A$ and $x/E = H \nsubseteq (x/R^e)^c \cap H = B^c$. Thus, $x \notin A^\downarrow$
and $x \notin B^{c \downarrow}$. This means that 
$x \notin A^\downarrow \cup  B^{c \downarrow}$.

As we noted, if $R\circ R^{-1} = R^e$, then $\mathrm{RS}$ is a Stone algebra. Moreover,
it is a Kleene--Stone algebra in which complemented and sharp elements coincide. By Proposition~\ref{prop:BZstar)},
$\mathrm{RS}$ forms a PBZ*-lattice such that $\Neg(A,B) = (A,B)^*$.
\end{proof}

An \emph{antiortholattice} is a PBZ*-lattice with the property that $0$ and $1$ are its only sharp elements.

\begin{corollary}\label{antiortho}
Let $R$ be a quasiorder or a tolerance induced by an irredundant covering of $U$. We can define an antiortholattice on $\mathrm{RS}$ 
if and only if $R^e = U \times U$. In that case, operation $\Neg$ is given by
\begin{equation} \label{eq:all_negation}
\Neg (A,B) = \left\{
\begin{array}{ll}
  (U,U) & \text{if $(A,B) = (\emptyset ,\emptyset)$,} \\
  (\emptyset ,\emptyset) & \text{otherwise.}
\end{array} 
\right. 
\end{equation}
\end{corollary}

\begin{proof}
If $R$ be a quasiorder or a tolerance induced by an irredundant covering of $U$, $\mathrm{RS}$ is a complete distributive lattice.
By Theorem~\ref{Thm:PBZ_characterization},  $(\mathrm{DM(RS)}, \vee, \wedge, {\sim},\Neg, (\emptyset,\emptyset),(U,U))$ is a
PBZ-lattice if and only if there exists an equivalence $E$ extending $R$ and $\Neg (A,B)$ equals $(B^{c \downarrow}, B^{c \downarrow})$
for all $(A,B) \in \mathrm{RS}$. 

Assume now that $\mathrm{RS}$ forms an antiortholattice. Then $(\emptyset ,\emptyset)$ and $(U,U)$
are the only sharp elements. In view of Proposition~\ref{Prop:cent1} and Corollary~\ref{Cor:cent}, they coincide with the exact elements of $\mathrm{RS}$.
According to Proposition \ref{Prop:equivalence_class},  $(A,A)$ is exact  if and only if $A$ is a union of $R^e$-classes.
This implies that there is just one $R^e$-class, which is $U$. Thus, $R^e = U \times U$.

Conversely, assume that $R^e = U \times U$. Then $R^e$ is the only equivalence extending $R$. 
By Theorem~\ref{Thm:PBZ_characterization}, the only way to define an PBZ-lattice is by setting
$\Neg (A,B) = (B^{c \downarrow},B^{c \downarrow})$, where $^\downarrow$ is defined in terms of $R^e$.
Now $\Neg(\emptyset,\emptyset) = (\emptyset^{c \downarrow},\emptyset^{c \downarrow}) = (U^{\downarrow},U^{\downarrow}) = (U,U)$ and
if $(A,B) \neq (\emptyset,\emptyset)$, then $B^c \subset U$, $B^{c \downarrow} = \emptyset$, and
$\Neg(A,B) = (\emptyset,\emptyset)$. Thus, $\Neg$ is defined as in \eqref{eq:all_negation}.

Finally, we prove that $\Neg ((A,B) \wedge  {\sim} (A,B)) \leq \Neg (A,B) \vee \Neg {\sim} (A,B)$ holds for all $(A,B)\in RS$. 
We have two possibilities: (i) $(A,B) \wedge {\sim}(A,B) \neq (\emptyset,\emptyset)$ or
(ii) $(A,B) \wedge {\sim}(A,B) = (\emptyset,\emptyset)$.
In case (i), $\Neg ((A,B) \wedge {\sim}(A,B)) = (\emptyset,\emptyset)$ and the condition holds trivially. In case (ii),
we have that $(A,B)$ is complemented. By Proposition~\ref{Prop:cent1} and Corollary~\ref{Cor:cent},
$(A,B)$ is exact, that is, $(A,B) = (B,B)$. By Proposition~\ref{Prop:equivalence_class}, $B$ is a union of $R^e$-classes.
There is only one $R^e$-class, which is $U$. Thus, the only exact sets are $(\emptyset,\emptyset)$ and $(U,U)$. 
If $(A,B) = (\emptyset,\emptyset)$, then $\Neg(A,B) = (U,U)$. If $(A,B) = (U,U)$, then $\Neg {\sim} (A,B) = (U,U)$.
Also in this case, the required condition holds.
\end{proof}

\section*{Conclusions}
\label{sec:Conc}

In the case of rough sets induced by arbitrary binary relations, we knew quite a
little about their structure. Practically only the results presented in \cite{Umadevi2015} about the completion $\mathrm{DM(RS)}$ were known. 
In this work, we have extended this knowledge in the case of a reflexive relation by showing that $\mathrm{DM(RS)}$
forms a paraorthomodular lattice. Also, the connections between exact rough sets, central, sharp and complemented elements of $\mathrm{DM(RS)}$ were studied. 

Since we have shown the one-to-one correspondence between atomistic complete Boolean 
sublattices of $\mathrm{DM(RS)}$ and PBZ-lattices on $\mathrm{DM(RS)}$, 
it is clear that PBZ-lattices can always be defined on $\mathrm{DM(RS)}$; there
always exists the Boolean sublattice $\{(\emptyset,\emptyset), (U,U)\}$.
This opens the possibility to study the properties of the complete lattice
$\mathrm{DM(RS)}$ for some particular reflexive relations different from the
known cases.
If $R$ is a quasiorder or a tolerance induced by an irredundant covering,
then $\mathrm{RS}$ itself is a completely distributive lattice and any PBZ-lattice
definable on it can be induced by an equivalence relation extending $R$. For
a quasiorder $R$, the only PBZ*-lattice which can be built on it is
induced by the least equivalence $R^{e}$ containing $R$. We can see that in
these cases, the obtained algebraic structures are not richer than those
described in \cite{JPR13} and \cite{JarRad17}. 

In case of an equivalence relation $E$, the rough sets defined by $E$ form a
PBZ*-lattice whose negation $\Neg$ coincides to the pseudocomplementation operation
in the lattice $\mathrm{RS}$, and this is the only PBZ*-lattice which can be
defined on it.

As we have pointed out, not all pseudo-Kleene
algebras are isomorphic to $\mathrm{DM(RS)}$ for some reflexive relation.
In the further studies, our aim is to find how pseudo-Kleene algebras $\mathrm{DM(RS)}$ can be characterized among all the others.

\section*{Acknowledgements}
We thank the anonymous referees for their valuable remarks on our manuscript.


\end{document}